\documentclass[a4paper,11pt]{article}
\usepackage{stmaryrd}
\usepackage{graphicx}
\usepackage{amsfonts}
\usepackage{amsmath,amsthm}
\usepackage{amssymb}
\usepackage{mathrsfs}
\usepackage{indentfirst}
\usepackage{titlesec}
\usepackage{caption2}
\usepackage{color}
\usepackage[normalem]{ulem}
\usepackage[english]{babel}
\usepackage{cite}
\usepackage{bbm}

\newtheorem{theorem}{Theorem}[section]
\newtheorem{definition}{Definition}[section]

\newtheorem{lemma}[theorem]{Lemma}
\newtheorem{corollary}[theorem]{Corollary}

\begin{document}

\title{On the fast convergence of random perturbations of the gradient flow}

\author{Jiaojiao Yang
\thanks{School of Mathematics and Statistics, Anhui Normal University, Wuhu, 241002, P.~R. China. Email: \texttt{y.jiaojiao1025@yahoo.com}. J.Y was supported by National Natural Science Foundation of China (Grant No.11801199), Natural Science Foundation of Anhui Province (CN)(1908085QA30) and Anhui Normal University Foundation(CN)(903-751819).} \ , \
Wenqing Hu
\thanks{Department of Mathematics and Statistics, Missouri University of Science and Technology
(formerly University of Missouri, Rolla). Rolla, MO, USA. Email: \texttt{huwen@mst.edu}. W.H was supported by a University of Missouri Research Board Grant.} \ , \
Chris Junchi Li
\thanks{University of California, Berkeley. Berkeley, CA, USA.
Email: \texttt{junchi.li.duke@gmail.com}} \
}

\date{}

\maketitle

\begin{abstract}
We consider in this work small random perturbations (of multiplicative noise type) of the gradient flow. We prove that under mild conditions, when the potential function is a Morse function with additional strong saddle condition, the perturbed gradient flow converges to the neighborhood of local minimizers in $O(\ln (\varepsilon^{-1}))$ time on the average, where $\varepsilon$ is the scale of the random perturbation. Under a change of time scale, this indicates that for the diffusion process that approximates the stochastic gradient method, it takes (up to logarithmic factor) only a linear time of inverse stepsize to evade from all saddle points. This can be regarded as a manifestation of fast convergence of the discrete-time stochastic gradient method, the latter being used heavily in modern statistical machine learning.
\end{abstract}



\section{Introduction}

\subsection{Setup and the main results}
Let $F:\mathbb{R}^n\to \mathbb{R}$ be a function in class $\mathbf{C}^{(3)}$ with bounded first and second derivatives and let
$\nabla F$ be its gradient vector field on $\mathbb{R}^n$. In this paper we consider small
random perturbations (of multiplicative noise type)
of the gradient flow associated with the function $F$.
Let $Y^\varepsilon_t$ be defined as the solution to
the stochastic differential equation (SDE) with small parameter $\varepsilon>0$:

\begin{equation}\label{Eq:SGDPerturbation}
dY^\varepsilon_t
=
- \nabla F(Y^\varepsilon_t)dt + \varepsilon \sigma(Y^\varepsilon_t)dW_t
\ , \ Y_0^\varepsilon=x \ .
\end{equation}
Here $\sigma(\bullet)$ is an $n\times n$ matrix-valued function with bounded coefficients
in class $\mathbf{C}^{(2)}$ and bounded first derivatives of these coefficients,
such that the diffusion matrix $a(x)=\sigma(x)\sigma^T(x)$
is uniformly positive definite (which indicates that the noise is sufficiently omnidirectional).

In the case where $\varepsilon = 0$, the SDE \eqref{Eq:SGDPerturbation} reduces to the gradient flow ODE:

\begin{equation}\label{Eq:GradientFlow}
\dfrac{d}{dt}S^t x=-\nabla F(S^t x) \ , \ S^0x=x \ .
\end{equation}

It is well-known, see e.g., ~Sec.~10.2 of \cite{[BORKAR]}, that the equilibria of the gradient flow differential equation \eqref{Eq:GradientFlow} are the critical points of $F(\bullet)$ (i.e., the gradient vector evaluated at the critical point is 0). Classical results state that under mild regularity conditions for $F(\bullet)$, every solution $S^t x$ of the gradient flow \eqref{Eq:GradientFlow} on $\mathbb{R}^n$ exists for all $t \ge 0$, and $S^t x$ converges to a connected component of the set of critical points of $F$ as $t \to \infty$. Additionally, following the trajectories of gradient flow \eqref{Eq:GradientFlow}
the function value $F$ is nonincreasing, but they may be trapped in the neighborhood of a saddle point for a substantial amount of time. To see the first conclusion, note that for any solution of ODE, $F(\bullet)$ itself serves as a Lyapunov function:
$$
\frac{d F(S^t x)}{dt}
=
\nabla F(S^t x)\cdot  \frac{d S^t x}{dt}
=
-| \nabla F(S^t x)|^2 \le 0
\ ,
$$
and hence $F (S^t x)$ is a non-increasing function.

Our goal in this paper is to use tools from stochastic analysis to prove the following result: adding a small amount of random perturbation (multiplicative noise type) enables fast evasion from saddle points and lands the process in a neighborhood of the set of local minimum points. Let $x^*$ be a local minimum point of $F(\bullet)$ in the sense that for some open
neighborhood $U(x^*)$ of $x^*$:
$$
x^*
=
\arg\min\limits_{x\in U(x^*)} F(x)
.
$$
Assume $F(x)>F(x^*)$,
then given a sufficiently small constant $0<e<F(x) - F(x^*)$ (the exact condition for $e$ is specified in \eqref{Eq:Condition-e}), we aim to characterize the distribution of the stopping time $T^\varepsilon_x$ defined as
\begin{equation}\label{Eq:HittingTimeSGDPerturbationManySaddle}
T^\varepsilon_x \equiv \inf\{t\geq 0: Y^\varepsilon_0=x \ , \ F(Y^\varepsilon_t ) \leq  F(x^*) + e\}
\end{equation}
as $\varepsilon \to 0^+$. Equivalently, when starting from $Y^\varepsilon_0=x$, $T^\varepsilon_x$ is the first time the perturbed gradient flow dynamics $Y^\varepsilon_t$ leads to at least the amount of $F(x) - F(x^*) - e$ decay in function value of $F$. In a hand-waiving manner, our main result (Theorem \ref{Thm:MajorResultSGDPerturbationConvergenceTime}) can be formulated as the following

\paragraph{Main Result.}
\textit{Let the function $F: \mathbb{R}^n\rightarrow \mathbb{R}$ be a Morse function satisfying the ``strong saddle condition'' that will be specified in Section 1.2. When $\varepsilon\to 0^+$, $\mathbf{E} T^\varepsilon_x$ is asymptotically bounded by
$C \ln (\varepsilon^{-1})$, where $C>0$ is some constant that is determined by the landscape of the function $F$.}

\

This result is significant in the sense that, when there exist many saddle points on the landscape of the function
$F(\bullet)$ the ODE \eqref{Eq:GradientFlow}
risks being trapped at around saddle points. Nevertheless as long as the function satisfy some additional landscape property,
by adding small random perturbation $\varepsilon \sigma(Y^\varepsilon_t) dW_t$ the random dynamics \eqref{Eq:SGDPerturbation}
pays a merely $\ln(\varepsilon^{-1})$ factor of time (multiplied by a constant) to
enter into a neighborhood of the set of local minimum points.

\subsection{Strong saddle condition}

To detail the landscape condition of $F(\bullet)$, we first remind the readers of the definition of Morse function, as the following
\begin{definition}\label{Def:Morse}
A function $F:\mathbb{R}^n\to \mathbb{R}$ is a {\it Morse function} if it is infinitely differentiable and has all its
critical points being {\it non-degenerate}, i.e., for each critical point $x$, all
eigenvalues of the Hessian matrix $\nabla^2 F(x)$ are nonzero.
\end{definition}

Morse functions admit a local quadratic re-parametrization at each critical point,
which is the content of the so-called {\it Morse Lemma} \cite[Lemma 2.2]{[Milnor1963]}.
To ensure that the perturbed gradient flow escapes from saddle points, we introduce the following ``strict saddle property"
(compare with \cite{[GeEtAl1503]}, \cite{[sun2015nonconvex]}) as follows.

\begin{definition}[strict saddle property]\label{Def:StrictSaddleProperty}
Given fixed $\gamma_1>0$ and $\gamma_2>0$, we say a Morse function $F$ defined on $\mathbb{R}^n$ satisfies the ``strict saddle property" if each point $x\in \mathbb{R}^n$ belongs to one of the following:
(i) $|\nabla F(x)|\geq \gamma_2 >0$ \ ;
(ii) $|\nabla F(x)|< \gamma_2$ and $\lambda_{\min}(\nabla^2 F(x))\leq -\gamma_1 <0$ \ ;
(iii) $|\nabla F(x)|< \gamma_2$ and $\lambda_{\min}(\nabla^2 F(x))\geq \gamma_1>0$.
\end{definition}

Here $\lambda_{\min}(\nabla^2 F(x))$ stands for the minimum eigenvalue of the Hessian matrix $\nabla^2 F(x)$. Notice that since $F$ is defined on the whole space $\mathbb{R}^n$, to satisfy the strict saddle property, we usually have $F(x)\rightarrow\infty$ or $F(x)\rightarrow -\infty$ as $|x|\rightarrow \infty$. We will call a saddle point $x\in \mathbb{R}^n$ of the function $F$ a ``strict saddle" if Definition \ref{Def:StrictSaddleProperty} (ii)
 holds at $x$. Thus a Morse function $F$
that satisfies the strict saddle property has all its saddle points being strict saddle points.

For the sake of proof, it is natural to assume that all eigenvalues of the Hessian $\nabla^2 F$ at critical points are uniformly bounded away from 0. This leads to our new notion of ``strong saddle property" as follows.

\begin{definition}[strong saddle property]\label{Def:StrongSaddleProperty}
Let the Morse function $F(\bullet)$ satisfy the strict saddle property with parameters
$\gamma_1>0$ and $\gamma_2>0$. We say the Morse function $F(\bullet)$ satisfy the
``strong saddle property" if
for some $\gamma_3 > 0$ and any $x \in \mathbb{R}^n$ such that
$\nabla F(x) = 0$,
all eigenvalues $\lambda_i$, $i=1,2,...,n$ of the Hessian $\nabla^2 F(x)$ at $x$
satisfying (ii) in Definition \ref{Def:StrictSaddleProperty} are bounded away from zero by some
$\gamma_3>0$ in absolute value,
i.e., $|\lambda_i|\geq \gamma_3>0$ for any $1\leq i \leq n$.
\end{definition}

We will call a saddle point $x\in \mathbb{R}^n$ of the function $F$ a ``strong saddle" if Definition \ref{Def:StrongSaddleProperty} holds at $x$. Thus a Morse function $F$ that satisfies the strong saddle property has all its saddle points that are strong saddle points. Throughout this paper we will work under Definition \ref{Def:StrongSaddleProperty} for the function $F$.

\subsection{Linerization of the gradient flow near a strong saddle point}

By the classical Hartman-Grobman Theorem (see \cite[\S 13]{[Arnold]}), for any strong saddle point $O$
that we consider, there exists an open neighborhood $U$ of $O$, and a $\mathbf{C}^{(0)}$ homeomorphism mapping $f: U\rightarrow \mathbb{R}^n$,
such that the gradient flow \eqref{Eq:GradientFlow} is mapped by $f$ into a linear flow. The homeomorphism $f$ is called
a (linear) conjugacy mapping. It turns out, that this mapping
can be taken to be $h$-H\"{o}lder
continuous for some $0<h\leq 1$ that depends only on the Hessian $\nabla^2 F(O)$ at $O$ (see \cite{[StrengthenedHartmanGrobman]}).

We put here an additional assumption regarding the conjugacy mapping.

\,

\textbf{Linerization Assumption.} \textit{The homeomorphism $f$ provided by the Hartman-Grobman Theorem
can be taken to be $\mathbf{C}^{(2)}$.}

\,

This assumption is needed in the proof of Lemma \ref{Appendix:Lm:ExitTimeExitDistributionPart2}, which is important for Theorem \ref{Thm:Kifer1981Strengthened}.
It is known that a sufficient condition for the
validity of the $\mathbf{C}^{(2)}$ Linerization Assumption is the so called non-resonance condition (see, for example,
the Sternberg linerization Theorem \cite[Theorem 6.6.6]{[KatokHasselblatt]}). We refer the reader to \cite[Theorem 6.6.6]{[KatokHasselblatt]}
for more discussions about our Linerization Assumption.

\subsection{Relation with statistical machine learning}

In this subsection we draw a connection between the perturbed gradient flow defined in \eqref{Eq:SGDPerturbation} and the stochastic gradient method in modern statistical machine learning. Stochastic Gradient Descent (SGD, \cite{BottouSGD2010}, \cite{Bach-MoulinesSGD2013NIPS}, \cite{Flammarion-BachCOLT2015}) with constant learning rate is a stochastic analogue of the gradient descent algorithm, aiming at finding the local or global minimizers or maximizers of the function expectation parameterized by some random variable. One can schematically formulate the optimization problem as follows: under some initial distribution, we target at finding a local minimum point $x^*$ of the expectation of function $F(x) \equiv \mathbf{E}[ F(x;\zeta) ]$, where the indexed random variable $\zeta$ follows some prescribed distribution. If the regularity conditions on the exchangeability of gradient operator and expectation operator hold, e.g., when $\zeta$ is supported on a finite set, the gradient method updates via the iteration
$$x_t
=
x_{t-1} - \beta
\mathbf{E}[ \nabla F(x_{t-1};\zeta) ]
.$$
However when the scale of the problem is extremely large, the access of expected gradients $\mathbf{E}[ \nabla F(z_{t-1};\zeta) ]$ are often expensive, and SGD prevails due to its one-query of noisy gradient at each iteration. In particular, SGD iteration often takes the following form
\begin{equation}\label{Eq:DiscreteSGD}
x_t
=
x_{t-1} - \beta
\nabla F(x_{t-1}; \zeta_t) \
,
\end{equation}
where $\beta$ is the fixed step-size (called the learning rate \footnote{In the classical Robbins-Monro algorithm (see \cite{Robbins-Monroe}), the coefficient $\beta$ may be taken as a decreasing sequence. However, in many modern machine learning practices, the coefficient $\beta$ is chosen to be constant.}) and $\{\zeta_t\}$ are i.i.d.~random variables that have the same distribution as $\zeta$. In particular, in the case of training a deep neural network,
the random variable $\zeta$ samples  mini-batches $\mathcal{B}$ with size $m$ ($m\leq n$) uniformly from an index set $\{1,2,...,n\}$:
$\mathcal{B}\subset \{1,2,...,n\}$ and $|\mathcal{B}|=m$. In this case, given loss functions $f_1(x)=\ell_{(a_1, b_1)}(x), ..., f_n(x)=\ell_{(a_n, b_n)}(x)$ on
training data $(a_1, b_1), ..., (a_n, b_n)$, we have $\nabla F(x; \zeta)=\dfrac{1}{m}\sum\limits_{i\in \mathcal{B}}\nabla f_i(x)$ and $\mathbf{E}[\nabla F(x; \zeta)]=\dfrac{1}{n}\sum\limits_{i=1}^n \nabla f_i(x)$. Due to its advantage in scalable data, SGD gains tremendous popularity in solving many large-scale statistical machine learning problems in the age of {\it Big Data}.

To analyze them, we associate the SGD with a diffusion process, i.e., the solution to a stochastic differential equation. Conversely, the SGD iteration can be viewed as a discrete-time, numerical scheme of such diffusion process. This point of view is closely related to stochastic approximation theory
and the weak convergence theory of Markov processes, and has been thoroughly investigated in many previous works, notably the classical monographs \cite{Kushner-Yin}, \cite{[BORKAR]},
 \cite{BenvenisteEtAlAdaptiveAlgorithms}. In particular, \cite[Chapters 7-9]{Kushner-Yin} describes the weak convergences of discrete iterations to SDE.
Here we follow the setting of the recent work by \cite{JunchiEtAlSDE2017} and analyze a constant small step-size, continuous-time analogue of the stochastic gradient method that takes the form
\begin{equation}\label{Eq:SGD}
dX_t
=
-\beta \nabla F(X_t)dt + \beta \sigma(X_t)dW_t
\ ,  \ X_0=x \ ,
\end{equation}
where $\sigma(z) = \left[\text{Var}\left(\nabla F(z;\zeta)\right)\right]^{1/2}$ is a positive semidefinite matrix. Equivalently, one can accelerate a factor of $\beta^{-1}$ and obtain for $Y^{\sqrt{\beta}}_t \equiv X_{t/\beta}$,
$$
dY^{\sqrt{\beta}}_t
=
- \nabla F(Y^{\sqrt{\beta}}_t)dt + \sqrt{\beta} \sigma(Y^{\sqrt{\beta}}_t)dW_t
\ , \ Y^{\sqrt{\beta}}_0=x \ .
$$
In other words, \eqref{Eq:SGDPerturbation} holds for $\varepsilon=\sqrt{\beta}$. Therefore under the small step-size regime $\beta\downarrow 0$, the stochastic process $Y^{\sqrt{\beta}}_t$ defined in \eqref{Eq:SGDPerturbation} can be viewed as a continuous-time analogue of the SGD iteration, which addresses the stochasticity emerging from SGD.
We further refer to \cite{WeinanEtAlSDE2017},  \cite{WeinanEtAlSDE2018}, \cite{JunchiEtAlSDE2017} for the proof of the convergence of discrete SGD \eqref{Eq:DiscreteSGD} to its continuous version \eqref{Eq:SGD}. We also notice \cite{[SGD-continuous-time]} that deals with a problem also related to \eqref{Eq:SGD}. The equation \eqref{Eq:SGD} has been used as a continuous model for the discrete SGD iteration \eqref{Eq:DiscreteSGD} in many earlier works (see for example \cite{ZhanxingSGDAnisotropic}). Moreover, recent empirical studies further convinced that \eqref{Eq:SGD} is indeed a good model for the limiting process of SGD at least for large-batch training. Actually, in a most recent study~\cite{simsekli2019tail}, authors
argued that the SGD noise is heavy-tailed and non-Gaussian, thus it might not be the best
tool to approximate SGD using SDE driven by Brownian motion. This point of view has been empirically examined in \cite{NonGaussianitySGD-Noise}, where the authors found that the conclusions of \cite{simsekli2019tail} are debatebale, and indeed for batch sizes $256$ and above, the distribution of the noise is best described as Gaussian at least in the early phases of training. These evidences further consolidate our use of \eqref{Eq:SGD} as a continuous model for SGD.

Thus from statistical machine learning point of view we are interested in the asymptotic as $\beta \downarrow 0$ of the hitting time
\begin{equation}\label{Eq:HittingTimeSGDManySaddle}
\tau^\beta = \inf\{t: F(X_t) \leq  F(x^*) + e\} \ .
\end{equation}
In plain words, how much time (upper bound) does it require \textit{for the function value to decay} by at least $F(X_0) - F(x^*) + e$ using stochastic gradient descent?

By using our main result stated in Section 1.1, we observe the phenomenon that when the the noise term $\beta\sigma(X_t)dW_t$ in \eqref{Eq:SGD} is sufficiently omnidirectional (i.e. the diffusion matrix $a(x)=\sigma(x)\sigma^T(x)$ is uniformly positive definite) it enables {\it fast} evasion from saddle points. Our Corollary \ref{Corollary:MajorResultSGDConvergenceTime} in this paper suggests a time complexity $\mathbf{E}\tau^\beta \lesssim O(\beta^{-1}\ln (\beta^{-1}))$ as $\beta \downarrow 0$ in discrete-time SGD (compare with \cite{[GeEtAl1503]}, which requires $O(\beta^{-2})$ iteration). In other words, the hitting time to a point, in terms of function value, close to a local minimizer takes time that only introduces an additional factor of $O(\ln (\beta^{-1}))$ (compare with \cite{[GeEtAl1503]}).

It is worth mentioning that the previous work by Pemantle \cite{pemantle1990nonconvergence} proves from a gradient flow viewpoint that SGD can avoid all nondegenerate saddle points, whenever the diffusion matrix $a(x)$ is uniformly positive definite or ``omnidirectional'' (see also a recent work by \cite{lee2016gradient} for the pure gradient descent case). However, neither of these works provides an analysis on the convergence rates. Within the scope of this work, we concentrate on the diffusion process associated with the discrete-time SGD, and we do not aim to quantify the approximation errors when bridging the discrete-time and continuous-time versions. We believe that using tools from numerical stochastic differential equations, such approximation error is small compared to the noise of dynamics and can be estimated rigorously, which potentially inspires another thread of interesting works.

It is also interesting to notice a recent work \cite{[Pathway-Energy-Landscape]} related to the problem we are considering here. There the authors
discuss algorithms that construct pathway maps on energy landscapes which are essentially connecting the saddle points over the landscape.

\

The rest of this paper is organized as follows. In Section 2 we quantify the decay of function value as well as the exit time in the neighborhood of one specific strong saddle point. We turn to analyze the general scenario in Section 3 where there are finitely many saddles that the process encounters along its trajectory. We provide the precise statement of theorem in the final part of this paper.

\section{Analysis of the function decay and exit time in a neighborhood of one specific ``strong saddle"}

In this section we analyze the decay in function value of $F$ and the exit time asymptotic in a neighborhood of one specific
``strong saddle" point $x_s\in \mathbb{R}^n$. These problems have been discussed thoroughly in the literature on random perturbations
 of dynamical systems (see for example \cite{[Kifer1981]},
\cite{[Bakhtin2008SPA]}, \cite{[BakhtinHeteroclinic]}, \cite{[MonterBakhtinNormalForm]}, among others), however, these results do not directly carry to our setting, as we may have a chain of saddle points that may occur on the landscape of the function $F$. It is worth mentioning here that the setting closest to our work is \cite{[BakhtinHeteroclinic]}, where the author is also working with a chain of saddle points. The difference is that there these saddle points are connected by a hetero-clinic network, while in our setting, the saddles can be connected by gradient dynamics of \eqref{Eq:GradientFlow}.

We consider the trajectory
of \eqref{Eq:SGDPerturbation} in a neighborhood of one specific strong saddle point $x_s$.
Without loss of generality we can assume that
the saddle point $x_s=0$ is the origin $O$ inside $U \subset \mathbb{R}^n$: a bounded connected open set with smooth boundary.
 Let $G\subset U$
be a domain with smooth boundary, so that $O$ is the only isolated saddle inside $G$.

Let $O$ belong to the interior of $G$ and there are no other equilibriums in $G$.
Let $A=\nabla^2 F(O)$. We assume that $-\lambda_1=-\lambda_2=...=-\lambda_{q}<-\lambda_{q+1}\leq...\leq -\lambda_{p}
<0<\lambda_{p+1}\leq ... \leq \lambda_n$
are the (real) eigenvalues of $A$, and $\lambda_1\geq \gamma_1>0$.

We shall work with the time-rescaled and perturbed process $Y_t^\varepsilon$ as in \eqref{Eq:SGDPerturbation}.
Let
\begin{equation}\label{Eq:HittingTimeSGDPerturbation}
\tau^\varepsilon_x=\inf\{t>0: Y_0=x \ , \ Y^\varepsilon_t\in \partial G\} \ .
\end{equation}
Let us also consider the corresponding
deterministic gradient flow \eqref{Eq:GradientFlow}.
We introduce a decomposition

$$G\cup \partial G=O\cup A_1\cup A_2\cup A_3 \ ,$$
where $A_1$ is a set of points $x\in G\cup \partial G$ such that if $x\in A_1$ then $S^ux\in G$
for $u>s$ and $S^u x\not \in G\cup \partial G$ if $u\leq s$ for some $s=s(x)\leq 0$ and $S^t x\rightarrow O$ as $t\rightarrow \infty$;
$A_2$ is a set of points $x\in G\cup \partial G$ such that if $x\in A_2$ then $S^u x\in G$ for $u<s$
and $S^u x\not\in G\cup \partial G$ if $u>s$ for some $s=s(x)\geq 0$ and $S^t x\rightarrow O$ as $t\rightarrow -\infty$;
$A_3$ is a set of points $x\in G\cup \partial G$ such that if $x\in A_3$ then $S^u x\in G$ provided $s_1<u<s_2$
and $S^u x\not\in G\cup \partial G$ if either $u>s_2$ or $u<s_1$ for some $s_1=s_1(x)\leq 0$ and $s_2=s_2(x)\geq 0$.
Note that $A_1$ and $A_3$ may be simultaneously empty.

If $x\in A_2\cup A_3$, then $S^t x$ leaves $G$ after some time, so that there is a finite
\begin{equation}\label{Eq:DeterministicFiniteExitTime}
t(x)=\inf\{t>0 : S^t x\in \partial G\} \ .
\end{equation}

Let us denote by $\Gamma_{\text{max}}$ the eigenspace of $A$ which corresponds to the eigenvalues
$-\lambda_1,...,-\lambda_q$. It is known that (see \cite[Chapter 9]{[HartmanODE]}), there exists a $q$-dimensional sub-manifold $W_{\text{max}}$
tangent to $\Gamma_{\text{max}}$ at $O$ and is invariant with respect to $S^t$. We see that $Q_{\text{max}}=W_{\text{max}}\cap \partial G$
is not empty. If $q>1$ then $Q_{\text{max}}$ is a sub-manifold of ($q-1$)-dimensions on the boundary $\partial G$. If
$q=1$ then $Q_{\text{max}}$ consists of two points.

Let $U\subset G$ be an open neighborhood of the saddle point $O$.
Let $U$ be chosen so small that
$\text{dist}(U\cup \partial U, \partial G)>0$. Let
$t(x)$ be defined as in \eqref{Eq:DeterministicFiniteExitTime}. Set
$$\partial G_{U\cup \partial U \rightarrow {\text{out}}}=\{S^{t(x)}x \text{ for some }
x\in U \cup \partial U \text{ with finite } t(x)\}\cup Q_{\text{max}} \ .$$
For small $\mu>0$ we let
$$Q^\mu=\{x\in \partial G \ , \ \text{dist}(x, \partial G_{U\cup \partial U\rightarrow \text{out}})<\mu\} \ .$$

Then we have

\begin{theorem}\label{Thm:Kifer1981Strengthened}
For any $r>0$, there exist some $\varepsilon_0>0$ so that for all $x\in U\cup \partial U$ and all $0<\varepsilon<\varepsilon_0$ we have
\begin{equation}\label{Eq:ExpectedExitTimeStrengthened}
\dfrac{\mathbf{E}_x\tau_x^\varepsilon}{\ln(\varepsilon^{-1})}\leq \dfrac{1}{\lambda_1}+r \ .
\end{equation}
For any small $\mu>0$ and any $\rho>0$, there exist some $\varepsilon_0>0$ so that for all $x\in U\cup \partial U$ and all $0<\varepsilon<\varepsilon_0$
we have
\begin{equation}\label{Eq:ExitDistributionStrenghtened}
\mathbf{P}_x(Y_{\tau_x^\varepsilon}^\varepsilon\in Q^\mu)\geq 1-\rho \ ,
\end{equation}
which indicates that with high probability, the exit position $Y_{\tau_x^\varepsilon}^\varepsilon$ is close to $\partial G_{U\cup \partial U\rightarrow \text{out}}$.
\end{theorem}

\textit{Proof.} Our strategy is to decompose $U\cup \partial U$ into disjoint regions according to their distances to $A_1\cup O$. We prove the result when the initial point $x$ is close to $A_1\cup O$ within distance $\varepsilon>0$ in Lemma \ref{Appendix:Lm:ExitTimeExitDistributionPart1}, where we choose $0<\varepsilon<\varepsilon_0$ and $\varepsilon_0=\min(\varepsilon_0(r), \varepsilon_0(\rho, \mu))$. Then we prove the result when the initial point $x$ is gradually farther away from $A_1\cup O$, but still is ``close" to $A_1\cup O$ when $\varepsilon\rightarrow 0$, in the sense that its distance to $A_1\cup O$ is between $\varepsilon^{1/2^k}$ and $\varepsilon^{1/2^{k-1}}$ for $k=1,2,...$. This is done in Lemma
\ref{Appendix:Lm:ExitTimeExitDistributionPart2}, where for each $k$, the statement is true for $0<\varepsilon<\varepsilon_0$ such that $\varepsilon_0=\min(\varepsilon_0(r, k), \varepsilon_0(\rho, \mu, k))$. Finally, we show that when $k$ is larger than some finite $k_0=k_0(r,\rho,\mu)$, if we start the process $Y_t^\varepsilon$ from an initial point $x$ that is of a distance at least $\varepsilon^{1/2^{k_0}}$ from $A_1\cup O$, the statement is still true, provided that $0<\varepsilon<\varepsilon_0$ and $\varepsilon_0=\min(\varepsilon_0(r, k_0), \varepsilon_0(\rho, \mu, k_0))$. This is proved in Lemma \ref{Appendix:Lm:ExitTimeExitDistributionPart3}. Thus we choose
$$0<\varepsilon_0<\min(\varepsilon_0(r), \varepsilon_0(\rho, \mu), \varepsilon_0(r, 1), \varepsilon_0(\rho, \mu, 1), ..., \varepsilon_0(r, k_0), \varepsilon_0(\rho, \mu, k_0)) \ ,$$
where $k_0=k_0(r,\rho,\mu)$, so that we can conclude the result of this Theorem.
$\square$

\

Let us denote

\begin{equation}\label{Appendix:Eq:FirstNeighborhoodStableManifold}
(A_1\cup O)_{[0,\varepsilon)}=\{x: \text{dist}(x,A_1\cup O)<\varepsilon\} \ ,
\end{equation}
and for $k=1,2,...$ we denote
\begin{equation}\label{Appendix:Eq:ArbitraryNeighborhoodStableManifold}
(A_1\cup O)_{[\varepsilon^{1/2^{k-1}}, \varepsilon^{1/2^k})}=\{x: \varepsilon^{1/2^{k-1}}\leq \text{dist}(x,A_1\cup O)<\varepsilon^{1/2^k}\} \ .
\end{equation}

We provide the following three Lemmas, which settle the proof of Theorem \ref{Thm:Kifer1981Strengthened}.

\begin{lemma}\label{Appendix:Lm:ExitTimeExitDistributionPart1}
(i) For any $r>0$, there exist some $\varepsilon_0=\varepsilon_0(r)$ such that for any $0<\varepsilon<\varepsilon_0$
and any $x\in U\cup \partial U$ with $x\in (A_1\cup O)_{[0, \varepsilon)}$ we have
$$\mathbf{E}\tau_x^\varepsilon\leq \left(\dfrac{1}{\lambda_1}+r\right)\ln(\varepsilon^{-1}) \ .$$

(ii) For any $\rho>0$ and $\mu>0$, there exist some $\varepsilon_0=\varepsilon_0(\rho, \mu)$ such that for any $0<\varepsilon<\varepsilon_0$
and any $x\in U\cup \partial U$ with $x\in (A_1\cup O)_{[0, \varepsilon)}$ we have
$$\mathbf{P}(Y^\varepsilon_{\tau_x^\varepsilon}\in Q^\mu)\geq 1-\rho \ .$$
\end{lemma}

\textit{Proof}. If $x\in A_1\cup O$, then the proof of this lemma is more or less the same as in \cite{[Kifer1981]}.
If $0<\text{dist}(x, A_1\cup O)<\varepsilon$, one can modify the proof in \cite{[BakhtinHeteroclinic]},
and the argument is more or less the same as in the proof
of the next Lemma, so that we refer to the proof of the next Lemma
for more details. $\square$

\begin{lemma}\label{Appendix:Lm:ExitTimeExitDistributionPart2}
(i) For any $k\in \mathbb{N}$ and any $r>0$, there exist some $\varepsilon_0=\varepsilon_0(r,k)$ such that for any $0<\varepsilon<\varepsilon_0$
and any $x\in U\cup \partial U$ with $x\in (A_1\cup O)_{[\varepsilon^{1/2^{k-1}}, \varepsilon^{1/2^k})}$ we have
$$\mathbf{E}\tau_x^\varepsilon\leq \left(\dfrac{1}{2^{k-1}\lambda_1}+r\right)\ln(\varepsilon^{-1}) \ .$$

(ii) For any $k\in \mathbb{N}$ and any $\rho>0$, any $\mu>0$, there exist some $\varepsilon_0=\varepsilon_0(\rho,\mu,k)$ such that for any $0<\varepsilon<\varepsilon_0$
and any $x\in U\cup \partial U$ with $x\in (A_1\cup O)_{[\varepsilon^{1/2^{k-1}}, \varepsilon^{1/2^k})}$ we have
$$\mathbf{P}(Y^\varepsilon_{\tau_x^\varepsilon}\in Q^\mu)\geq 1-\rho \ .$$
\end{lemma}

\textit{Proof}. Let $V\supset U$ be an open neighborhood of $O$.
Let $f: V\rightarrow \mathbb{R}^n$ be the $\mathbf{C}^{(2)}$ conjugacy that maps the flow $S^t x$ in \eqref{Eq:GradientFlow}
within $V$ into linear dynamics $\dot{y}=\Lambda y$, $f(O)=0\in \mathbb{R}^n$. Explicitly, the solution of the linear flow
$\dot{y}=\Lambda y$, $y(0)=y_0=(y_0^1,...,y_0^n)$ is given by $y(t)=(y^1(t),...,y^n(t))$ such that $y^l(t)=y_0^l e^{\lambda_l t}$ for $l=1,2,...,n$
\footnote{In Section 2 we use the eigenvalues $\lambda_1,...,\lambda_p, -\lambda_{p+1},...,-\lambda_n$. But
here for simplicity we denote the eigenvalues by $\lambda_1,...,\lambda_n$ and we understand that $\lambda_1,...,\lambda_p$
are positive eigenvalues and $\lambda_{p+1},...,\lambda_n$ are negative eigenvalues of $\Lambda$.}.
We will denote the linear flow corresponding to $S^t$
to be $S^t_\Lambda$, so that $f(S^t x)=S^t_\Lambda (f(x))$.

We set $\mathcal{Y}_\varepsilon(t)=f(Y^\varepsilon_t)$. According to the argument that leads to (8.1)
at the beginning of Section 8 of \cite{[BakhtinHeteroclinic]}, we have

\begin{equation}\label{Appendix:Eq:BakhtinHeteroclinic8-1}
d\mathcal{Y}_\varepsilon(t)=\Lambda \mathcal{Y}_\varepsilon(t)+\varepsilon B(\mathcal{Y}_\varepsilon(t))dW(t)+\varepsilon^2 C(\mathcal{Y}_\varepsilon(t))dt \ ,
\end{equation}
where $B$ and $C$ are continuous and bounded in a neighborhood $f(V)=\{f(Y); Y\in V\}$ of the origin $0$, and $B$
is non-degenerate.

From the assumption that $x\in (A_1\cup O)_{[\varepsilon^{1/2^{k-1}}, \varepsilon^{1/2^k})}$ we know that the initial condition
\begin{equation}\label{Appendix:Eq:BakhtinHeteroclinic8-2}
\mathcal{Y}_\varepsilon(0)=y_0+\varepsilon^{1/2^k}\xi_\varepsilon \ ,
\end{equation}
where $y_0\in f(A_1\cup O)$ and $\xi_\varepsilon$ is a point on $f(V)$ such that $\varepsilon^{1/2^k}\leq \text{dist}(\xi_\varepsilon, f(A_1\cup O))<1$.

Let $\widetilde{V}\subset V$ and $W=f(V)$, $\widetilde{W}=f(\widetilde{V})$, so that
$$\widetilde{W}=\{y\in \mathbb{R}^n: |y^j|\leq R, j=1,2,...,n\}\subset W$$
for some small $R>0$. Thus $y_0^l=0$ for all $l\leq p$.

We would like to study the system \eqref{Appendix:Eq:BakhtinHeteroclinic8-1} with initial data given by \eqref{Appendix:Eq:BakhtinHeteroclinic8-2}.
Let us define the stopping times
$$\sigma_{W,\varepsilon}=\inf\{t\geq 0: \mathcal{Y}_\varepsilon(t)\in \partial W\} \ ,$$
$$\sigma_\varepsilon=\inf\{t\geq 0: \max(|\mathcal{Y}_{\varepsilon}^1(t)|,...,|\mathcal{Y}_{\varepsilon}^q(t)|)\geq R\}\wedge \sigma_{W, \varepsilon} \ .$$
It turns out, and will be clear from the proof that follows, that
$\mathbf{P}(\sigma_\varepsilon<\sigma_{W,\varepsilon})\rightarrow 1 \text{ as } \varepsilon\rightarrow 0$.

We have the solution $\mathcal{Y}_\varepsilon(t)$ in mild form
$$\mathcal{Y}_\varepsilon(t)=e^{\Lambda t}\mathcal{Y}_\varepsilon(0)+\varepsilon e^{\Lambda t}\int_0^t e^{-\Lambda s}B(\mathcal{Y}_\varepsilon(s))dW_s+\varepsilon^2 e^{\Lambda t}\int_0^t e^{-\Lambda s}C(\mathcal{Y}_\varepsilon(s))ds \ .$$

Fix some number $a\in \left(0,\dfrac{1}{2^k}\right)$. For every $l=1,2,...,n$ let us define
$$\tau_\varepsilon^l=\inf\{t: |\mathcal{Y}_\varepsilon^l(t)-(S^t_\Lambda y_0)^l|=\varepsilon^{a}\} \ ,$$
and we define
$$\tau_\varepsilon=\min\{\tau_\varepsilon^l: l=1,2,...,n\} \ .$$

In very much the same way as the derivation of Lemmas 8.2-8.4 (see Section 11 of \cite{[BakhtinHeteroclinic]})
we see that we have

\begin{equation}\label{Appendix:Eq:BakhtinHeteroclinicLemma8-2}
\sup\limits_{t\geq 0}|\mathcal{Y}_\varepsilon(\tau_\varepsilon\wedge t)-S^{t\wedge \tau_\varepsilon}_\Lambda y_0|\stackrel{\mathbf{P}}\rightarrow 0 \text{ as } \varepsilon\rightarrow 0 \ ,
\end{equation}

\begin{equation}\label{Appendix:Eq:BakhtinHeteroclinicLemma8-3}
\tau_\varepsilon\stackrel{\mathbf{P}}\rightarrow \infty \text{ as } \varepsilon\rightarrow 0 \ ,
\end{equation}

\begin{equation}\label{Appendix:Eq:BakhtinHeteroclinicLemma8-4}
\mathbf{P}(\tau_\varepsilon=\min(\tau_\varepsilon^1,...,\tau_\varepsilon^q))\rightarrow 1 \text{ as } \varepsilon\rightarrow 0 \ .
\end{equation}

By using the same reasoning as in \cite{[BakhtinHeteroclinic]}, which derives equation (8.9)
in \cite{[BakhtinHeteroclinic]}, we see that we have, for $l\leq p$,
\begin{equation}\label{Appendix:Eq:BakhtinHeteroclinic8-9}
\mathcal{Y}_\varepsilon^l(\tau_\varepsilon)=e^{\lambda_l \tau_\varepsilon}(\varepsilon^{1/2^k}\xi_\varepsilon^k+\varepsilon N_\varepsilon^k(\tau_\varepsilon)+o_\mathbf{P}(\varepsilon)) \ .
\end{equation}
Here $N_\varepsilon(\tau_\varepsilon)=(N^1_\varepsilon(\tau_\varepsilon),...,N^n_\varepsilon(\tau_\varepsilon))$ and $N_\varepsilon(\tau_\varepsilon)\stackrel{\text{Law}}\rightarrow N_0$
as $\varepsilon\rightarrow 0$, where $N_0$ is a Gaussian random vector defined in (8.3) of \cite{[BakhtinHeteroclinic]}. We have used
the notation that $\phi(\varepsilon)=o_{\mathbf{P}}(\psi(\varepsilon))$ for any families of random variables
$\phi(\varepsilon)$, $\psi(\varepsilon)$, $\varepsilon>0$ such that $\frac{\phi(\varepsilon)}{\psi(\varepsilon)}\stackrel{\mathbf{P}}\rightarrow 0$ as $\varepsilon\rightarrow 0$.

Put $\mathcal{K}_\varepsilon^l=\mathcal{Y}_\varepsilon^l(\tau_\varepsilon)e^{-\lambda_l \tau_\varepsilon}\varepsilon^{-1/2^k}=\xi_\varepsilon^l+\varepsilon^{1-1/2^k}N_\varepsilon^l(\tau_\varepsilon)
+o_\mathbf{P}(\varepsilon^{1-1/2^k})$, we rewite \eqref{Appendix:Eq:BakhtinHeteroclinic8-9} as

\begin{equation}\label{Appendix:Eq:BakhtinHeteroclinic8-10}
\mathcal{Y}_\varepsilon^l(\tau_\varepsilon)=e^{\lambda_l\tau_\varepsilon}\varepsilon^{1/2^k}\mathcal{K}_\varepsilon^l \ .
\end{equation}

From \eqref{Appendix:Eq:BakhtinHeteroclinicLemma8-4} we know that

\begin{equation}\label{Appendix:Eq:BakhtinHeteroclinic8-11}
\mathbf{P}(\max(|\mathcal{Y}^1_\varepsilon(\tau_\varepsilon)|,...,|\mathcal{Y}^q_\varepsilon(\tau_\varepsilon)|)=\varepsilon^a)\rightarrow 1 \ .
\end{equation}
So, we have $e^{\lambda_1\tau_\varepsilon}\varepsilon^{1/2^k}\max(|\mathcal{K}_\varepsilon^1|,...,|\mathcal{K}_\varepsilon^q|)=\varepsilon^a$, which gives, with probability
approaching $1$, as $\varepsilon\rightarrow 0$,

\begin{equation}\label{Appendix:Eq:BakhtinHeteroclinic8-12}
\tau_\varepsilon-\frac{\frac{1}{2^k}-a}{\lambda_1}\ln(\varepsilon^{-1})=-\dfrac{\ln (\max(|\mathcal{K}_\varepsilon^1|,...,|\mathcal{K}_\varepsilon^q|))}{\lambda_1} \ .
\end{equation}

Since we know $\varepsilon^{1/2^k} \leq \text{dist}(\xi_\varepsilon, f(A_1\cup O))<1$, we see that
for all $x\in (A_1\cup O)_{[\varepsilon^{1/2^{k-1}}, \varepsilon^{1/2^k})}$, the corresponding $\xi_\varepsilon=(\xi_\varepsilon^1,...,\xi_\varepsilon^n)$
satisfies $\varepsilon^{1/2^k}\leq |\xi_\varepsilon^l|<1$ for $l=1,2,...,n$.

For $k\geq 2$ we have, as $\varepsilon>0$ is small, with probability approaching $1$,
$$
\begin{array}{ll}
|\mathcal{K}_\varepsilon^l| & =|\xi_\varepsilon^l+\varepsilon^{1-1/2^k}N_\varepsilon^l(\tau_\varepsilon)
+o_\mathbf{P}(\varepsilon^{1-1/2^k})|
\\
&=\varepsilon^{1/2^k}\left|\varepsilon^{-1/2^k}\xi_\varepsilon^l+\varepsilon^{1-1/2^k-1/2^k}N_\varepsilon^l(\tau_\varepsilon)+o_{\mathbf{P}}(\varepsilon^{1-1/2^{k}-1/2^k})\right|
\\
& \stackrel{(a)}{\geq} \varepsilon^{1/2^k}\left|1-\varepsilon^{1-1/2^{k-1}}|N_\varepsilon^l(\tau_\varepsilon)|-o_\mathbf{P}(\varepsilon^{1-1/2^{k-1}})\right| ,
\end{array}
$$
so that
\begin{equation}\label{Appendix:Eq:BoundForcKTermsAuxiliarykgeq2}
-\ln|\mathcal{K}_\varepsilon^l|\leq \dfrac{1}{2^k}\ln(\varepsilon^{-1})+o_\mathbf{P}(1) \ .
\end{equation}
Here in (a) we used the fact that $\varepsilon^{-1/2^k}|\xi_{\varepsilon}^l|\geq 1$, $1-1/2^k-1/2^k=1-1/2^{k-1}$, $|a+b|\geq \left||a|-|b|\right|$, and with high probability $\varepsilon^{1-1/2^{k-1}}|N_\varepsilon^l(\tau_\varepsilon)|$ is smaller than $1$ as $\varepsilon\rightarrow 0$, since
$N_\varepsilon(\tau_\varepsilon)\stackrel{\text{Law}}{\rightarrow} N_0$ for a Gaussian random vector $N_0$ defined in (8.3) of [4] (with fixed covariances).

For $k=1$ we have, as $\varepsilon>0$ is small,
$$
\begin{array}{ll}
|\mathcal{K}_\varepsilon^l| & =|\xi_\varepsilon^l+\varepsilon^{1/2}N_\varepsilon^l(\tau_\varepsilon)
+o_\mathbf{P}(\varepsilon^{1/2})|
\\
& \geq \varepsilon^{1/2}\left||N_\varepsilon^l(\tau_\varepsilon)|-1-o_\mathbf{P}(\varepsilon^{1/2})\right| ,
\end{array}
$$
so that
\begin{equation}\label{Appendix:Eq:BoundForcKTermsAuxiliarykeq1}
-\ln|\mathcal{K}_\varepsilon^l|\leq \dfrac{1}{2}\ln(\varepsilon^{-1})+
\left|\ln(\left||N_\varepsilon^l(\tau_\varepsilon)|-1-o_\mathbf{P}(\varepsilon^{1/2})\right|)\right| \ .
\end{equation}

Equations \eqref{Appendix:Eq:BoundForcKTermsAuxiliarykeq1} and \eqref{Appendix:Eq:BoundForcKTermsAuxiliarykgeq2}
 give that for any $r>0$, there exist some $\varepsilon_0>0$ uniformly for all $\xi_\varepsilon$
with $\varepsilon^{1/2^k} \leq \text{dist}(\xi_\varepsilon, f(A_1\cup O))<1$,
as $0<\varepsilon<\varepsilon_0$, with probability approaching $1$,

\begin{equation}\label{Appendix:Eq:BoundForcKTerms}
\dfrac{-\ln(\max(|\mathcal{K}_\varepsilon^1|, ...,|\mathcal{K}_\varepsilon^q|))}{\ln(\varepsilon^{-1})} \leq \dfrac{1}{2^k}+r \ .
\end{equation}

Combining \eqref{Appendix:Eq:BakhtinHeteroclinic8-12} and \eqref{Appendix:Eq:BoundForcKTerms} we can conclude that,
for any $r>0$, there exist some
$\varepsilon_0>0$ uniformly for all $\xi_\varepsilon$ with $\varepsilon^{1/2^k} \leq \text{dist}(\xi_\varepsilon, f(A_1\cup O))<1$,
as $0<\varepsilon<\varepsilon_0$,

\begin{equation}\label{Appendix:Eq:BoundExpectedTravelTimeAroundSaddle}
\dfrac{\mathbf{E}\tau_\varepsilon}{\ln(\varepsilon^{-1})} \leq \dfrac{\frac{1}{2^{k-1}}-a}{\lambda_1}+r \ .
\end{equation}

Moreover, in very much the same way as the derivation of Lemmas 8.5 and 8.7 of \cite{[BakhtinHeteroclinic]}, we know that
for $q+1\leq l\leq p$ we have
\begin{equation}\label{Appendix:Eq:BakhtinHeteroclinicLemma8-5}
\mathcal{Y}_\varepsilon^l(\tau_\varepsilon)=\varepsilon^{-\frac{\lambda_l}{\lambda_1}(\frac{1}{2^k}-a)+\frac{1}{2^k}}
[\max(|\mathcal{K}_\varepsilon^1,...,\mathcal{K}_\varepsilon^q|)]^{-\frac{\lambda_l}{\lambda_1}}\mathcal{K}_\varepsilon^l \ ,
\end{equation}
and for $p+1\leq l\leq n$ and any $\beta>0$ we have
\begin{equation}\label{Appendix:Eq:BakhtinHeteroclinicLemma8-7}
\mathcal{Y}_\varepsilon^l(\tau_\varepsilon)=\varepsilon^{-\frac{\lambda_l}{\lambda_1}(\frac{1}{2^k}-a)}[\max(|\mathcal{K}_\varepsilon^1,...,\mathcal{K}_\varepsilon^q|)]^{-\frac{\lambda_l}{\lambda_1}}
(y_0^l+\varepsilon^{\frac{1}{2^k}}\xi_\varepsilon^l)+o_\mathbf{P}(\varepsilon^{1-\beta}) \ .
\end{equation}

Let us now consider the process $\mathcal{Y}_\varepsilon(t)$ after $\tau_\varepsilon$. Let $\bar{\mathcal{Y}}_\varepsilon(t)=\mathcal{Y}_\varepsilon(\tau_\varepsilon+t)$.
Then we have the mild form of the solution $\bar{\mathcal{Y}}_\varepsilon(t)$:

\begin{equation}\label{Appendix:Eq:BakhtinHeteroclinic8-15}
\bar{\mathcal{Y}}_\varepsilon(t)=e^{\Lambda t}\left(\mathcal{Y}_\varepsilon(\tau_\varepsilon)+\varepsilon\int_0^t e^{-\Lambda s}B(\bar{\mathcal{Y}}_\varepsilon(s))dW_s
+\varepsilon^2\int_0^t e^{-\Lambda s}C(\bar{\mathcal{Y}}_\varepsilon(s))ds\right) \ .
\end{equation}

Set $\bar{\tau}_\varepsilon=\inf\{t: \max(|\bar{\mathcal{Y}}_\varepsilon^1(t)|,...,|\bar{\mathcal{Y}}_\varepsilon^q(t)|)=R\}$, then we have, in parallel to (8.16)
of \cite{[BakhtinHeteroclinic]}, that

\begin{equation}\label{Appendix:Eq:BakhtinHeteroclinic8-16}
\max(|\bar{\mathcal{Y}}_\varepsilon^1(\bar{\tau}_\varepsilon)|,...,|\bar{\mathcal{Y}}_\varepsilon^q(\bar{\tau}_\varepsilon)|)
=e^{\lambda_1\bar{\tau}_\varepsilon}\varepsilon^{a}(1+\eta_\varepsilon) \ ,
\end{equation}
with $\eta_\varepsilon\stackrel{\mathbf{P}}\rightarrow 0$ as $\varepsilon\rightarrow 0$. Thus \eqref{Appendix:Eq:BakhtinHeteroclinic8-16} implies that

\begin{equation}\label{Appendix:Eq:BakhtinHeteroclinic8-17}
\bar{\tau}_\varepsilon=\dfrac{a}{\lambda_1}\ln(\varepsilon^{-1})+\dfrac{1}{\lambda_1}\ln\left(\dfrac{R}{1+\eta_\varepsilon}\right) \ .
\end{equation}

Back to part (i) of this Lemma that we are proving, we see that $\sigma_\varepsilon=\tau_\varepsilon+\bar{\tau}_\varepsilon$
and $0\leq \mathbf{E}\tau_x^\varepsilon-\mathbf{E} \sigma_\varepsilon\leq C$ for some constant $C>0$ independent of $\varepsilon$ when $\varepsilon$ is small.
This, together with \eqref{Appendix:Eq:BoundExpectedTravelTimeAroundSaddle}, \eqref{Appendix:Eq:BakhtinHeteroclinic8-17},
prove part (i) of this Lemma.

Regarding part (ii) of this Lemma, it is readily checked from \eqref{Appendix:Eq:BakhtinHeteroclinicLemma8-5} and
\eqref{Appendix:Eq:BakhtinHeteroclinicLemma8-7}, in parallel to (8.18) of \cite{[BakhtinHeteroclinic]}, that as $\varepsilon\rightarrow 0$,

\begin{equation}\label{Appendix:Eq:BakhtinHeteroclinic8-18}
\sup\limits_{t\leq \bar{\tau}_\varepsilon}|\bar{\mathcal{Y}}_\varepsilon(t)-S^t_\Lambda (\varepsilon^a v)|\stackrel{\mathbf{P}}\rightarrow 0
\end{equation}
where $v\in \text{span}(v_1,...,v_q)$ with $v_1,...,v_q$ being normal unit vectors in the eigenvector directions corresponding to
$\lambda_1,...,\lambda_q$. This, together with the fact that as $\varepsilon\rightarrow 0$,
\begin{equation}\label{Appendix:Eq:BakhtinHeteroclinicLemma9-1}
\sup\limits_{t\in [0,T]}|\mathcal{Y}_\varepsilon(t)-S^t_\Lambda y_0|\stackrel{\mathbf{P}}\rightarrow 0 \ ,
\end{equation}
show part (ii) of this Lemma. $\square$

\begin{lemma}\label{Appendix:Lm:ExitTimeExitDistributionPart3}
(i) For any $r>0$, there exist some $k_0=k_0(r)\in \mathbb{N}$ such that, there exist some
$\varepsilon_0=\varepsilon_0(r,k_0)$ such that for any $0<\varepsilon<\varepsilon_0$
and any $x\in U\cup \partial U$ with $\text{dist}(x,A_1)\geq \varepsilon^{1/2^{k_0}}$ we have
$$\mathbf{E}\tau_x^\varepsilon\leq 2r\ln(\varepsilon^{-1}) \ .$$

(ii) For any $\rho>0$, any $\mu>0$, there exist some $k_0=k_0(\rho,\mu)\in \mathbb{N}$ such that,
there exist some $\varepsilon_0=\varepsilon_0(\rho,\mu,k_0)$ such that for any $0<\varepsilon<\varepsilon_0$
and any $x\in U\cup \partial U$ with $\text{dist}(x,A_1)\geq \varepsilon^{1/2^{k_0}}$ we have
$$\mathbf{P}(Y^\varepsilon_{\tau_x^\varepsilon}\in Q^\mu)\geq 1-\rho \ .$$
\end{lemma}

\textit{Proof}. For some $k_0\in \mathbb{N}$ let us consider the domain
$D=\{x\in G: \text{dist}(x, (A_1\cup O))> \varepsilon^{1/2^{k_0-1}}\}$. For any initial point $x\in D$, the deterministic flow
$S^t x$ hits $\partial G$ in a penetrating manner (see Definition \ref{Def:RegularANDPenetratingExit}) within time
$t(x,\partial G)\leq \dfrac{C}{2^{k_0}}\ln(\varepsilon^{-1})$, where $C>0$ is a constant independent of $\varepsilon$. From here, by using
an $\varepsilon$-dependent version of the arguments in Lemma \ref{Lm:FiniteTimeComparisonGradientFlowANDSGDPerturbation}, we conclude part (i)
of this Lemma. Part (ii) of this lemma follows the same way as \eqref{Eq:FiniteTimeComparisonGradientFlowANDSGDPerturbationExitPosition}
in Lemma \ref{Lm:FiniteTimeComparisonGradientFlowANDSGDPerturbation}.
$\square$

\

\section{The general case: finitely many saddles}

Let us consider the case when the
Morse function $F$ has finitely many strong saddles $O_1,...,O_k$ (according to Definition \ref{Def:StrongSaddleProperty}).
Since the problem formulated in Section 1 is only about
function value decay so that the process $Y_t^\varepsilon$ hits a local minimum point of $F$, we can assume,
without loss of generality, that the saddles $O_1,...,O_k$ are ordered in such a way that
$F(O_1)>F(O_2)>...>F(O_k)$, and the point $x^*$ is a local minimum point of $F$ with $F(x^*)<F(O_k)$.
We will work within the basin $U(x^*)$ of $x^*$, so that $O_1,...,O_k$ are the only critical points
besides $x^*$ in $U(x^*)$.
Let us work with the perturbed process $Y_t^\varepsilon$ in \eqref{Eq:SGDPerturbation}. Let $Y_0^\varepsilon=x\in U(x^*)$
and $F(x)>F(O_1)$. Given a fixed small number $e>0$, such that
\begin{equation}\label{Eq:Condition-e}
F(O_k)>F(x^*)+e \ ,
\end{equation}
let
$$T^\varepsilon_x=\inf\{t\geq 0: Y_0^\varepsilon=x \ , \  F(Y_t^\varepsilon)\leq F(x^*)+e\} \ .$$
Our goal is to give an asymptotic estimate of $T^\varepsilon_x$.

For each saddle point $O_i$, we consider a nested pair of open neighborhoods $U_i\subsetneq V_i$ containing $O_i$.
Let us first pick all $V_i$, $i=1,2,...,k$ in such a way that
$O_i\in V_i$ is the only stationary point inside $V_i$, and
for any $i\neq j$, we have $V_i\cap V_j=\emptyset$. Let us then pick $U_i\subset V_i$
such that $\text{dist}(\partial U_i, \partial V_i)>0$ for all $i=1,2,...,k$. For each saddle point $O_i$,
let us denote by $\Gamma_{i, \text{max}}$ the eigenspace of $\nabla^2 F(O_i)$ which corresponds to the negative
eigenvalues of $\nabla^2 F(O_i)$ with largest absolute value.
For each open neighborhood $V_i$,
as in Section 2, we can construct a submanifold $W_{i, \text{max}}$ that is tangent
to $\Gamma_{i, \text{max}}$ at $O_i$ and is invariant with respect to $S^t$.
Let $Q_{i, \text{max}}=W_{i, \text{max}}\cap \partial V_i$. By the classification of points in $V_i\cup \partial V_i$
as $O_i$, $A_1$, $A_2$,
$A_3$ as in Section 2, we see that for any point $x\in  U_i\cup \partial U_i$,
either there exist some finite $t(x)$ such that $S^{t(x)}x\in \partial V_i$, or $S^t x\rightarrow O_i$ as $t\rightarrow\infty$.
Let
$$\partial V_{i,  U_i\cup \partial U_i \rightarrow {\text{out}}}=\{S^{t(x)}x \text{ for some }
x\in U_i \cup \partial U_i \text{ with finite } t(x)\}\cup Q_{i, \text{max}} \ .$$

For any fixed positive number $h>0$, we introduce the level curves
$$C_i(F(O_i)\pm h)=\{x\in \mathbb{R}^n: F(x)=F(O_i)\pm h\} \ .$$

Given $h>0$ sufficiently small, let us introduce the open domain
$$D((F(O_i)+h, F(O_{i-1})-h))=\{x\in \mathbb{R}^n: F(O_i)+h<F(x)<F(O_{i-1})-h\} \ .$$
If $i=1$, we can symbolically let $F(O_0)=\infty$, so that
$$D((F(O_1)+h, F(O_0)-h))=D((F(O_1)+h, \infty))=\{x\in \mathbb{R}^n: F(O_1)+h<F(x)\} \ .$$

Note that the curve $C_i(F(O_i)\pm h)$ and the domain $D((F(O_i)+h, F(O_{i-1})-h))$
may have several disconnected components, and only some of these components
have points that are close to $O_i$ ($O_{i-1}$). Also note that $C_i(F(O_i)-h)$ ($C_{i-1}(F(O_{i-1})+h)$) are the
lower (upper) boundaries of $D((F(O_i)+h, F(O_{i-1})-h))$.

We will now present a geometric lemma, that describes the geometry (see Figure 1) of the gradient flow
in a neighborhood of a non-degenerate strong saddle.

\begin{lemma}\label{Lm:GeometricLemma}
(Geometric lemma) Uniformly for all $i=1,2,...,k$, we have the following.

(i) Besides the basic assumptions of the nested neighborhoods $U_i$, $V_i$ that we imposed before,
we can pick $U_i$, $V_i$ in such a way that there exist some constant $\mu>0$ independent of $i$, and some constant
$h>0$ that may depend on the choice of $U_i$ and $V_i$, such that
for any $x\in \partial V_i$ with $\text{dist}(x, \partial V_{i, U_i\cup \partial U_i\rightarrow \text{out}})<\mu$, we have
$F(x)<F(O_i)-h$;

(ii) One can pick the constant $h>0$ so small that $C_i(F(O_i)+\frac{h}{2})\cap U_i\neq \emptyset$;

(iii) For fixed $h>0$ as in (ii), there exist some
$\kappa>0$ such that $|\nabla F(x)|\geq \kappa$ for any point $x\in D(F(O_i)+\frac{h}{2}, F(O_{i-1})-\frac{h}{2})$;

(iv) One can pick the constant $\kappa>0$ so small that for any point $x\in C_i(F(O_i)+\frac{h}{2})$,
and $x\not \in U_i$, there exist a tubular
neighborhood $\mathcal{T}(x)$, such that $\partial \mathcal{T}(x)\cap C_i(F(O_i)-h)\neq \emptyset$,
$\partial \mathcal{T}(x)\cap C_i(F(O_i)+h)\neq \emptyset$, and all other boundaries of $\mathcal{T}(x)$ are parallel to flow lines of \eqref{Eq:GradientFlow}.
Moreover, for any $y\in \mathcal{T}(x)$, we have $|\nabla F(y)|\geq \kappa$.

\end{lemma}

\textit{Proof.} By the classical Hartman-Grobman Theorem (\cite[\S 13]{[Arnold]}), within a sufficiently small neighborhood
of the strong saddle point $O_i$, the flow $S^t x$ in \eqref{Eq:GradientFlow} is $\mathbf{C}^{(0)}$-conjugate
to the linear flow with matrix $-\nabla^2 F(O_i)$. Since the claims stated in this Lemma are all invariant
under $\mathbf{C}^{(0)}$ deformation, we just have to prove this lemma under the framework of linear dynamics (see Figure 1).

(i) By the so-called Inclination Lemma \cite[Theorem 5.7.2]{[Brin-Stuck]} (also see Lemma 5.7.1 in the same reference),
one can pick $U_i$ small enough so that all flows starting from $U_i\cup \partial U_i$, if not attracted to $O_i$,
will approach the unstable manifold of the saddle point $O_i$, before they hit $\partial V_i$. In fact, the Inclination Lemma
indicates that the flow $S^t x$ will map any small disk that is intersecting transversally
with the stable manifold, with the dimension equal to that of the unstable manifold, to be close to the unstable manifold in
$\mathbf{C}^{(1)}$-norm.
Moreover, this disk is stretched in the unstable directions while it is being attracted to the unstable manifold.
Due to non-degeneracy of the saddle point $O_i$, by picking $V_i$ relatively large with respect to $U_i$
 we see that the statement follows;

(ii) Note that if $h>0$ satisfy (i), then any $\widetilde{h}>0$ such that $\widetilde{h}<h$ also satisfy (i).
Thus by picking $\widetilde{h}$ sufficiently small we have (ii);

(iii) Fix $h>0$ as in (ii), since one can find a neighborhood $N_i$ of $O_i$ inside $U_i$, such that $C_i(F(O_i)+\frac{h}{2})$ does not touch
this neighborhood, and $O_i$'s are isolated zeros of $|\nabla F(x)|$, we see that we have (iii);

(iv) We can take the neighborhood $N_i$ of $O_i$ inside $U_i$ to be so small that for any point
$x\in C_i(F(O_i)+\frac{h}{2})$ and $x\not\in U_i$, the flow line of \eqref{Eq:GradientFlow}
passing through $x$ and connecting $C_i(F(O_i)+h)$ to $C_i(F(O_i)-h)$, does not touch $N_i$. In fact,
since by (ii) we know that $C_i(F(O_i)+\frac{h}{2})\cap U_i\neq \emptyset$, we have, that for any $x\not\in U_i$
and $x\in C_i(F(O_i)+\frac{h}{2})$, the point $x$ lies outside of some invariant stable cone $K_\delta^s$ (see Figure 1).
The Inclination Lemma again implies that the forward flow $S^t x$ stretches $x$ along the unstable direction, thus keeps
the resulting flow away from $N_i$.
Thus the flow line of \eqref{Eq:GradientFlow}
passing through $x$ and connecting $C_i(F(O_i)+h)$ to $C_i(F(O_i)-h)$, does not touch $N_i$.
Further, the family of flow lines in a tubular neighborhood $\mathcal{T}(x)$ of this flow line
are connecting $C_i(F(O_i)+h)$ to $C_i(F(O_i)-h)$, and they do not touch $N_i$. For any $y$ on
this $\mathcal{T}(x)$, we will have $|\nabla F(y)|\geq \kappa>0$.
$\square$

%
%
%

\begin{figure}
\centering
\includegraphics[height=9cm, width=13.5cm]{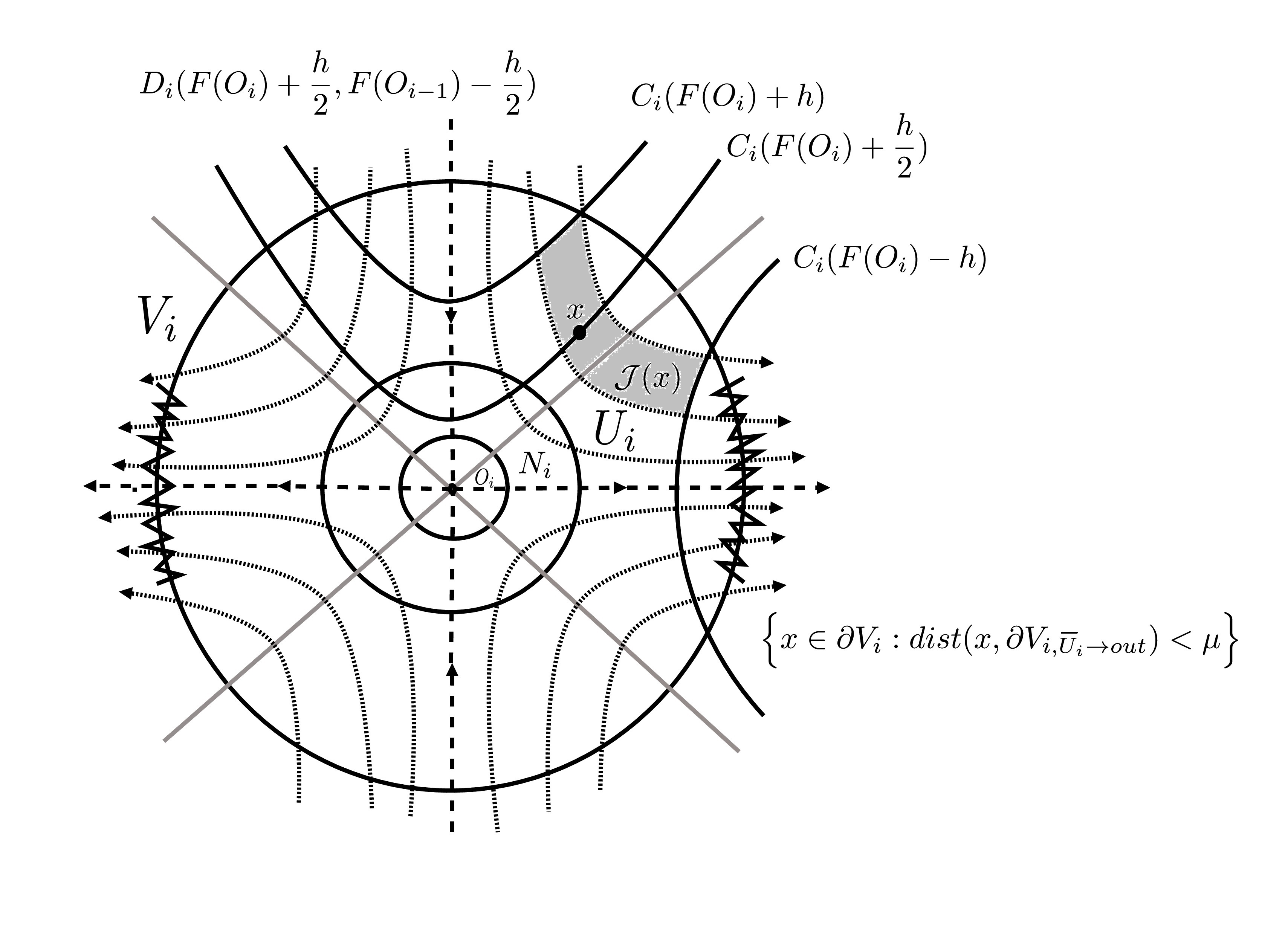}
\caption{Geometry of gradient flow dynamics near a saddle point.}
\end{figure}

\

Note that as $h\rightarrow 0$, we have $\kappa\rightarrow 0$, and it might happen that $\dfrac{h}{\kappa^2}\rightarrow\infty$.
However, in what follows we will pick some fixed $h>0$ and $\kappa>0$ as well as the neighborhoods of the saddle points
as in Lemma \ref{Lm:GeometricLemma}, and we let $\varepsilon\rightarrow 0$.

\

Let us provide now an auxiliary lemma, which is essentially adapted from \cite[Chapter 2, Lemma 3.1]{[FWbook]},
that will be used frequently in our subsequent arguments. Before we state the lemma, let us first introduce
some definition regarding exit behavior of the flow $\{S^t x: t \geq 0\}$ in \eqref{Eq:GradientFlow}.

Let $D\subset \mathbb{R}^n$ be an open domain in $\mathbb{R}^n$.
Let $x\in D$ and we consider the deterministic trajectory $\{S^t x: t\geq 0\}$.

\begin{definition}\label{Def:RegularANDPenetratingExit}
We say that the trajectory $\{S^t x: t\geq 0\}$ starting from initial condition $x\in D$,
exits $D$ in a ``regular manner" if
$$t(x, \partial D)=\inf\{t\geq 0: S^t x\not\in D\}<\infty \ , $$
and for some sufficiently small $\delta>0$ we have
\begin{equation}\label{Eq:RegularExit}
S^{t(x,\partial D)+\delta}(x)\not \in D\cup \partial D \ .
\end{equation}

We say that the trajectories $\{S^t x: t\geq 0, x\in D\}$ exit $D$
in a ``penetrating manner" if for any $x\in D$,
the trajectory $\{S^t x: t\geq 0\}$ exits $D$ in a regular manner,
such that there exist some constant $T_0>0$
with $t(x,\partial D)\leq T_0$ for all $x\in D$, and there exist some constant $c>0$,
such that for some sufficiently small $\delta>0$, we have
\begin{equation}\label{Eq:PenetratingEntrance}
\max\limits_{t(x,\partial D)\leq t\leq t(x,\partial D)+\delta}\text{dist} \ (S^t x, D\cup \partial D)\geq c >0 \ ,
\end{equation}
uniformly for all $x\in D$.

We say that the part of the boundary $\Gamma\subseteq \partial D$ is the ``exit piece" of $\{S^t x, x\in D, t\geq 0\}$
if $\{S^t x: x\in D, t\geq 0\}$ exit $D$ in a penetrating manner, and for any point $y\in \Gamma$, there exist some
$x\in D$ and $t(x,\partial D)\leq T_0$ such that $S^{t(x,\partial D)} x=y$.

\end{definition}

\begin{lemma}\label{Lm:FiniteTimeComparisonGradientFlowANDSGDPerturbation}
Let the domain $D\subset \mathbb{R}^n$ and initial point $x\in D$ be stated as in Definition \ref{Def:RegularANDPenetratingExit}.
Let $$t(x,\partial D)=\inf\{t\geq 0: S^t x\not \in D\} \ .$$
Let $$\tau^\varepsilon_x(\partial D)=\inf\{t\geq 0: Y_0=x \ , \ Y_t^\varepsilon\in \partial D\} \ .$$
If $t(x,\partial D)$ is finite for some choice of $x$ and $D$, and $\{S^t x: t\geq 0\}$ exits $D$ in a regular manner, then
$\tau_x^\varepsilon(\partial D) \rightarrow t(x,\partial D)$
in probability as $\varepsilon\downarrow 0$, i.e., for any given $\delta>0$,
\begin{equation}\label{Eq:FiniteTimeComparisonGradientFlowANDSGDPerturbationConvergenceInProbability}
\lim\limits_{\varepsilon\rightarrow 0}\mathbf{P}_x(|\tau^\varepsilon_x(\partial D)-t(x,\partial D)|>\delta)=0 \ .
\end{equation}

Moreover, if $\{S^t x, t\geq 0, x\in D\}$ exits $D$ in a penetrating manner, and for some small $r>0$ we have $x\in D^r\equiv \{x\in D, \text{dist}(x, \partial D)\geq r\}$,
then there exist some $\varepsilon_0>0$ that may depend on $D$ and $c$, and
some constant $C$ that may depend on $T_0$, but is independent of $\varepsilon_0$,
such that for any $0<\varepsilon<\varepsilon_0$ we have
\begin{equation}\label{Eq:FiniteTimeComparisonGradientFlowANDSGDPerturbationFiniteExpectedTime}
\sup\limits_{x\in D^r}\mathbf{E}_x \tau^\varepsilon_x(\partial D)\leq C \ .
\end{equation}

Finally, suppose $\{S^t x, t\geq 0, x\in D\}$ exits $D$ in a penetrating manner, and $\Gamma\subseteq \partial D$
is the exit piece as in Definition \ref{Def:RegularANDPenetratingExit}. Then for any open set $Q$
of $\partial D$ and $Q\supset\Gamma$, for any $\rho>0$ sufficiently small, there exist some $\varepsilon_0>0$ such that for any $0<\varepsilon<\varepsilon_0$, we have
\begin{equation}\label{Eq:FiniteTimeComparisonGradientFlowANDSGDPerturbationExitPosition}
\mathbf{P}_x(Y^\varepsilon_{\tau_x^\varepsilon(\partial D)}\in Q)\geq 1-\rho \ ,
\end{equation}
for all $x\in D^r$.

\end{lemma}

\textit{Proof.} If $\{S^t x, t\geq 0\}$ exits $D$ in a regular manner, then for every sufficiently small $\delta>0$
we have
$$S^{t(x,\partial D)-\delta}x\in D \ , \ S^{t(x,\partial D)+\delta}x\not\in D\cup \partial D \ .$$
Let $\delta_1>0$ be the distance of the trajectory segment $\{S^s x \ , \ s\in [0, t(x,\partial D)-\delta]\}$ from $\partial D$, and
$\delta_2>0$ be the distance of $S^{t(x,\partial D)+\delta}x$ from $\partial D$, and let $\bar{\delta}=\min(\delta_1,\delta_2)>0$. By
\cite[Chapter 2 , Theorem 1.2]{[FWbook]} we see that
$$\lim\limits_{\varepsilon\rightarrow 0}\mathbf{P}_x\left\{\sup\limits_{0\leq s \leq t(x,\partial D)+\delta}|Y^\varepsilon_s(x)-S^s x|>\bar{\delta}\right\}=0 \ .$$
Here $Y^\varepsilon_t(x)$ denotes the process $Y_t^\varepsilon$ in \eqref{Eq:SGDPerturbation} with $Y^\varepsilon_0=x$. This implies that
\begin{equation}\label{Eq:PointwiseEstimateHittingTime}
\lim\limits_{\varepsilon\rightarrow 0}\mathbf{P}_x(\tau^\varepsilon_x(\partial D)\in [t(x,\partial D)-\delta, t(x,\partial D)+\delta])=1 \ ,
\end{equation}
which is \eqref{Eq:FiniteTimeComparisonGradientFlowANDSGDPerturbationConvergenceInProbability}, and from here we see
$\tau^\varepsilon_x(\partial D)\rightarrow t(x,\partial D)$ in probability as $\varepsilon\downarrow 0$.

Now let us assume that $\{S^t x \ , \ t\geq 0 \ , \ x\in D\}$ exit $D$ in a penetrating manner. Let
$t(x,\partial D)\leq T_0$ for all $x\in D$ and $\max\limits_{t(x, \partial D)\leq t\leq t(x,\partial D)+\delta}\text{dist}
\ (S^t x, D\cup \partial D)\geq c >0$
uniformly for all $x\in D$ and some $\delta>0$. This implies that for every $\delta>0$, there exist some $\varepsilon_0>0$
such that for any $0<\varepsilon<\varepsilon_0$ we have
\begin{equation}\label{Eq:UniformEstimateHittingTime}
\mathbf{P}_x(|\tau^\varepsilon_x(\partial D)-t(x,\partial D)|>\delta)\leq \delta
\end{equation}
for all $x\in D^r$. Therefore
$$\sup\limits_{x\in D^r} \mathbf{P}_x(\tau^\varepsilon_x(\partial D)>2T_0)\leq \delta \ .$$

By strong Markov property of the process $Y_t^\varepsilon$ we see that

$$\begin{array}{ll}
\sup\limits_{x\in D^r} \mathbf{P}_x(\tau_x^\varepsilon(\partial D)>n\cdot 2 T_0) & = \sup\limits_{x\in D^r}
[\mathbf{P}_x(\tau^\varepsilon_{Y^\varepsilon_{(n-1)2T_0}}(\partial D)>2T_0|Y^\varepsilon_{(n-1)2T_0}\in D)\mathbf{P}_x(\tau_x^\varepsilon(\partial D)>(n-1)2T_0)]
\\
& \leq \delta \cdot \sup\limits_{x\in D^r}\mathbf{P}_x(\tau_x^\varepsilon(\partial D)>(n-1)2T_0) \ .
\end{array}$$

Therefore for every integer $n\geq 0$ we have, for every $x\in D^r$,

$$\mathbf{P}_x(\tau_x^\varepsilon(\partial D)> n \cdot 2T_0)\leq \delta^n \ .$$

Therefore we have

$$\begin{array}{ll}
\mathbf{E}_x\tau_x^\varepsilon(\partial D) & \leq \sum\limits_{n=0}^\infty (n+1)\cdot 2T_0\cdot \mathbf{P}_x(n\cdot 2T_0< \tau_x^{\varepsilon}(\partial D)\leq (n+1)\cdot 2T_0)
\\
& \leq \sum\limits_{n=0}^\infty (n+1)\cdot 2T_0\cdot \mathbf{P}_x(\tau_x^{\varepsilon}(\partial D)>n\cdot 2T_0)
\\
& \leq \sum\limits_{n=0}^\infty (n+1)\cdot 2T_0 \cdot \delta^n
\\
& \leq C
\end{array}$$
since we can choose $\delta>0$ to be sufficiently small. This implies \eqref{Eq:FiniteTimeComparisonGradientFlowANDSGDPerturbationFiniteExpectedTime}.

Finally, let us prove \eqref{Eq:FiniteTimeComparisonGradientFlowANDSGDPerturbationExitPosition}.
For any $\delta>0$ and any $\mu>0$, we have
$$\begin{array}{ll}
&\mathbf{P}_x\left(|Y^\varepsilon_{\tau_x^\varepsilon(\partial D)}-S^{t(x,\partial D)}x|>\delta\right)
\\
=  & \mathbf{P}_x\left(|Y^\varepsilon_{\tau_x^\varepsilon(\partial D)}-S^{t(x,\partial D)}x|>\delta \ , \ |\tau_x^\varepsilon(\partial D)-t(x,\partial D)|>\mu\right)
\\
& \ \ \ \ \
+\mathbf{P}_x\left(|Y^\varepsilon_{\tau_x^\varepsilon(\partial D)}-S^{t(x,\partial D)}x|>\delta \ , \ |\tau_x^\varepsilon(\partial D)-t(x,\partial D)|\leq \mu\right)
\\
\leq & \mathbf{P}_x\left(|\tau_x^\varepsilon(\partial D)-t(x,\partial D)|>\mu\right)
\\
& \ \ \ \ \
+\mathbf{P}_x\left(|Y^\varepsilon_{\tau_x^\varepsilon(\partial D)}-S^{t(x,\partial D)}x|>\delta \ , \ |\tau_x^\varepsilon(\partial D)-t(x,\partial D)|\leq \mu\right)
\\
\leq & \mathbf{P}_x\left(|\tau_x^\varepsilon(\partial D)-t(x,\partial D)|>\mu\right)
\\
& \ \ \ \ \
+\mathbf{P}_x\left(\sup\limits_{t(x,\partial D)-\mu \leq t\leq t(x,\partial D)+\mu}|Y^\varepsilon_t-S^{t(x,\partial D)}x|>\delta\right)
\\
\leq & \mathbf{P}_x\left(|\tau_x^\varepsilon(\partial D)-t(x,\partial D)|>\mu\right)
\\
& \ \ \ \ \
+\mathbf{P}_x\left(\sup\limits_{t(x,\partial D)-\mu \leq t\leq t(x,\partial D)+\mu}(|Y^\varepsilon_t-S^t x|+|S^t x-S^{t(x,\partial D)}x|)>\delta\right)
\\
\leq & \mathbf{P}_x\left(|\tau_x^\varepsilon(\partial D)-t(x,\partial D)|>\mu\right)
\\
& \ \ \ \ \
+\mathbf{P}_x\left(\sup\limits_{t(x,\partial D)-\mu \leq t\leq t(x,\partial D)+\mu}|Y^\varepsilon_t-S^t x|+
\sup\limits_{t(x,\partial D)-\mu \leq t\leq t(x,\partial D)+\mu}|S^t x-S^{t(x,\partial D)}x|>\delta\right) \ .
\end{array}$$

Due to the continuity of the flow $S^t x$, one can pick $\mu>0$ small enough such that for all $x\in D$ we have
$\sup\limits_{t(x,\partial D)-\mu \leq t\leq t(x,\partial D)+\mu}|S^t x-S^{t(x,\partial D)}x|<\dfrac{\delta}{2}$. Thus we have

$$\begin{array}{ll}
&\mathbf{P}_x\left(|Y^\varepsilon_{\tau_x^\varepsilon(\partial D)}-S^{t(x,\partial D)}x|>\delta\right)
\\
\leq & \mathbf{P}_x\left(|\tau_x^\varepsilon(\partial D)-t(x,\partial D)|>\mu\right)
+\mathbf{P}_x\left(\sup\limits_{t(x,\partial D)-\mu \leq t\leq t(x,\partial D)+\mu}|Y^\varepsilon_t-S^t x|>\dfrac{\delta}{2}\right) \ .
\end{array}$$
By \eqref{Eq:UniformEstimateHittingTime} as well as the fact that $t(x,\partial D)\leq T_0$ for all $x\in D^r$, we see that
the above estimate implies \eqref{Eq:FiniteTimeComparisonGradientFlowANDSGDPerturbationExitPosition}.
$\square$

\

Let us define

$$C(\pm h)=\bigcup\limits_{i=1}^k C_i(F(O_i)\pm h) \ .$$

For simplicity of notations let us also define $C_i(\pm h)=C_i(F(O_i)\pm h)$.

Let us define a sequence of stopping times (see Figure 2, one can also find similar constructions
in \cite[Chapter 6]{[FWbook]}, also see \cite{[Koralov2004]})

\begin{equation}\label{Eq:SequenceStoppingTime}
0=\sigma_0\leq \tau_1\leq \sigma_1\leq \tau_2\leq \sigma_2\leq ...
\end{equation}
such that
\begin{equation}\label{Eq:SequenceStoppingTimeTau}
\tau_j=\inf\{t>\sigma_{j-1}: Y_t^\varepsilon\in C(h/2)\} \ ,
\end{equation}

\begin{equation}\label{Eq:SequenceStoppingTimeSigma}
\sigma_j=\inf\{t>\tau_j: Y_t^\varepsilon \in C(-h)\} \ .
\end{equation}

%
%
%

\begin{figure}[h]
\centering
\includegraphics[height=12cm, width=12cm]{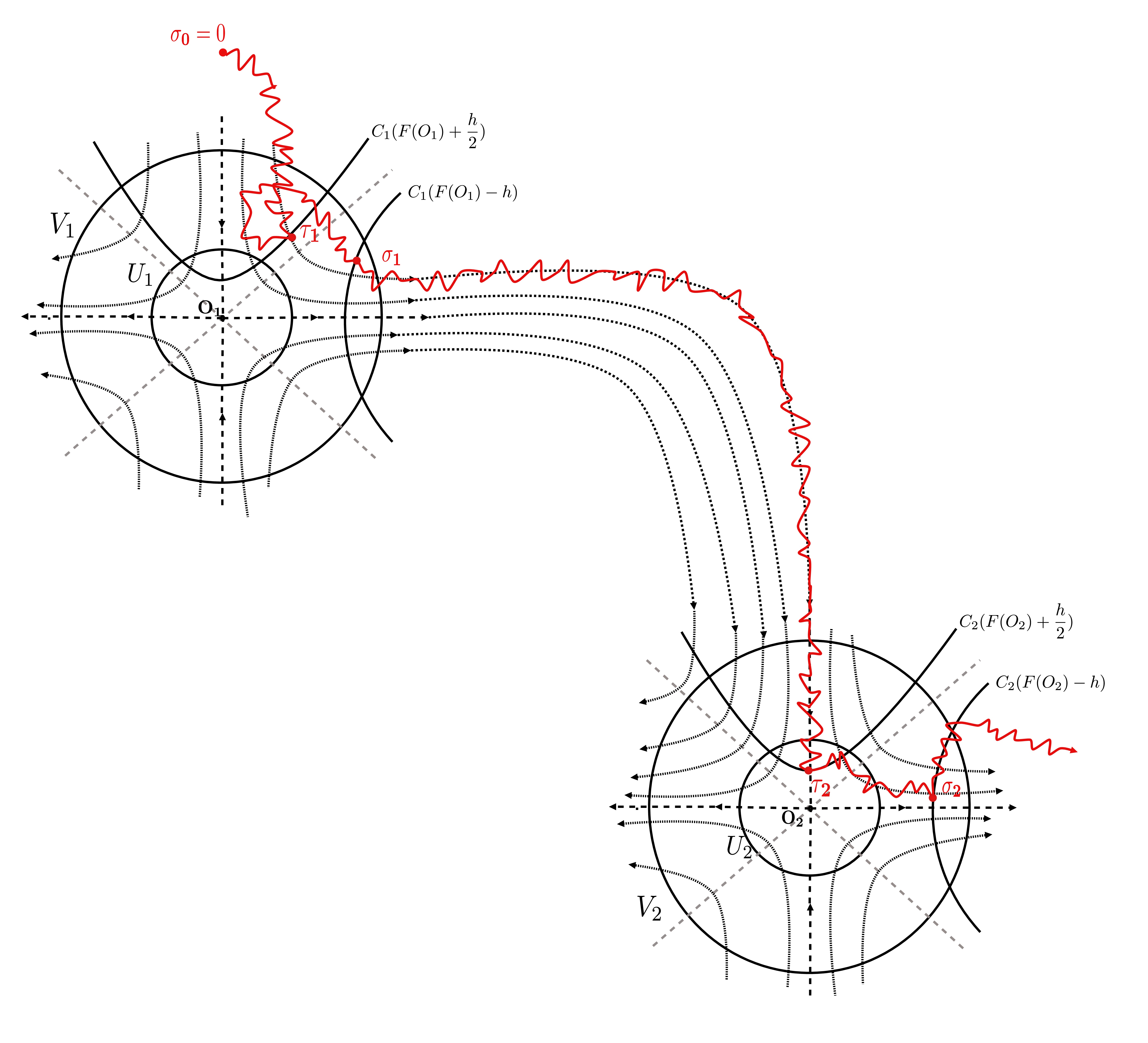}
\caption{Sequence of stopping times.}
\end{figure}

Let us start the process $Y_t^\varepsilon$ from $Y_0^\varepsilon=x\in \mathbb{R}^n$, such that $F(x)>F(O_1)$, and at the same time
$x\not\in V_1\cup \partial V_1$. Let us also assume that we have some $H>0$ such that for some small $r>0$ we have $H>F(x)+r$, and for all points $y\in \partial V_1$,
we have $H>F(y)+r$.

From our geometric lemma Lemma \ref{Lm:GeometricLemma}, we see that for any $x_1$ such that
$F(x_1)\geq F(O_1)+h/2$, we have
$|\nabla F(x_1)|\geq \kappa>0$. By the fact that $\dfrac{d}{dt}F(S^t x)=-|\nabla F(S^t x)|^2\leq -\kappa^2$, this implies that
\begin{equation}\label{Eq:HittingTimeUniformBoundTowardsBottomGradientFlow}
t(x,C_1(h/2))\leq \dfrac{F(x)-F(O_1)-h/2}{\kappa^2}<\dfrac{H-r-F(O_1)-h/2}{\kappa^2} \ .
\end{equation}
Moreover, since $S^t x$ is a gradient flow, the flow lines of $S^t x$ will be perpendicular to
 $C_1(h/2)$ when $x\in C_1(h/2)$. Due to the fact that $|\nabla F(x_1)|\geq \kappa>0$ when $x_1\in C_1(h/2)$, we see that
the flow $S^t x$ will hit $C_1(h/2)$ in a penetrating manner. Therefore
by Lemma \ref{Lm:FiniteTimeComparisonGradientFlowANDSGDPerturbation} we see that for such initial point $x$
we have the following.

\begin{lemma}\label{Lm:FirstTravelTimeFinite}
Given a small $r>0$, there exists some $\varepsilon_0>0$ uniformly for all $x$ with $H-r> F(x)\geq F(O_1)+h/2$,
such that for all $0<\varepsilon<\varepsilon_0$, there exist some finite $C>0$ independent of $\varepsilon$ such that
\begin{equation}\label{Eq:TravelTimeFinite}
\mathbf{E}_x \tau_1\leq C \ .
\end{equation}
\end{lemma}

Now we are ready to state and prove our main Theorem.

\begin{theorem}\label{Thm:MajorResultSGDPerturbationConvergenceTime}
Consider the process $Y_t^\varepsilon$ defined as in \eqref{Eq:SGDPerturbation}. Suppose the
Morse function $F(x)$
have $k$-strong saddle points (according to Definition
\ref{Def:StrongSaddleProperty}) $O_1$, ..., $O_k$ inside the basin $U(x^*)$, such that
$F(O_1)>...>F(O_k)$. Let $\gamma_1>0$
be chosen as in Definition \ref{Def:StrictSaddleProperty}. Let the
initial condition $Y_0=x\in U(x^*)$
be such that $F(x)>F(O_1)$. Let $x^*$ be the unique local minimum of
$F$ within $U(x^*)$ such that $F(x^*)<F(O_k)$. Then we have

(i) For any small $\rho>0$, with probability at least $1-\rho$, the process $Y_t^\varepsilon$ in
\eqref{Eq:SGDPerturbation} converges
to the minimizer $x^*$ for sufficiently small $\varepsilon$ after passing through all $k$ saddle points $O_1$, ..., $O_k$;

(ii) Consider the stopping time $T_x^\varepsilon$ defined in
\eqref{Eq:HittingTimeSGDPerturbationManySaddle}. Then as $\varepsilon\downarrow 0$, conditioned on the above convergence of $Y_t^\varepsilon$ to $x^*$,
we have
\begin{equation}\label{Eq:TotalExitTimeAsymptoticSGDPerturbationManySaddle}
\lim\limits_{\varepsilon \rightarrow 0}\dfrac{\mathbf{E} T_x^\varepsilon}{\ln (\varepsilon^{-1})}\leq
\dfrac{k}{\gamma_1} \ .
\end{equation}
\end{theorem}

\textit{Proof.}
By the strong Markov property of the process $Y_t^\varepsilon$, we see that
the process $\widetilde{Y}^\varepsilon_s=Y^\varepsilon_{\tau_1+s}$ can be viewed as an independent copy of the process
$Y_s^\varepsilon$ starting from $Y_0^\varepsilon=Y_{\tau_1}^\varepsilon$. By part (ii) of our Lemma \ref{Lm:GeometricLemma},
we can consider two cases, that either $Y^\varepsilon_{\tau_1}\in U_1\cup \partial U_1$ or $Y^\varepsilon_{\tau_1}\not\in U_1\cup \partial U_1$
but $Y^\varepsilon_{\tau_1}\in C_1(h/2)$.

Let us consider first the case when $Y^\varepsilon_{\tau_1}\in U_1\cup \partial U_1$.
Let $\mu>0$ be chosen as in Lemma \ref{Lm:GeometricLemma}.
Let
$$Q_1^\mu=\{x\in \partial V_1 \ , \ \text{dist}(x, \partial V_{1, U_1\cup \partial U_1\rightarrow \text{out}})<\mu\} \ .$$

Let $\widetilde{\sigma}_1=\inf\{t\geq \tau_1: Y_t^\varepsilon\in \partial V_1\}$.
Then by Theorem \ref{Thm:Kifer1981Strengthened} equation \eqref{Eq:ExitDistributionStrenghtened} we see that
 $$\lim\limits_{\varepsilon\rightarrow 0}\mathbf{P}_x(Y_{\widetilde{\sigma}_1}^\varepsilon\in Q_1^\mu |  Y_{\tau_1}^\varepsilon\in U_1\cup \partial U_1)=1 \ .$$

From the above we see that for any small $q>0$ there exist some $\varepsilon_0>0$ such that
for any $0<\varepsilon<\varepsilon_0$ we have
$$\mathbf{P}_x(Y_{\widetilde{\sigma}_1}^\varepsilon\in Q_1^\mu|  Y_{\tau_1}^\varepsilon\in U_1\cup \partial U_1)\geq 1-q \ .$$
By part (c) of our Lemma \ref{Lm:GeometricLemma} we know that for all $x\in Q_1^\mu$, we have $F(x)\leq F(O_1)-h$.
This implies that
$\{Y_{\widetilde{\sigma}_1}^\varepsilon \in Q_1^\mu\}\subseteq \{\widetilde{\sigma}_1\geq \sigma_1\}$. Therefore we have
$$\mathbf{P}_x(\widetilde{\sigma}_1\geq \sigma_1 |  Y_{\tau_1}^\varepsilon\in U_1\cup \partial U_1)\geq 1-q \ .$$
This indicates
\begin{equation}\label{Thm:MajorResultSGDPerturbationConvergenceTime:Proof:SigmatildeLargerSigma}
\lim\limits_{\varepsilon\rightarrow 0}\mathbf{P}_x(\widetilde{\sigma}_1\geq \sigma_1 |  Y_{\tau_1}^\varepsilon\in U_1\cup \partial U_1)=1 \ .
\end{equation}

From Theorem \ref{Thm:Kifer1981Strengthened} equation \eqref{Eq:ExpectedExitTimeStrengthened} we see that
\begin{equation}\label{Thm:MajorResultSGDPerturbationConvergenceTime:Proof:SigmatildeMinusTau}
\lim\limits_{\varepsilon\rightarrow 0}
\dfrac{\mathbf{E}_x (\widetilde{\sigma}_1-\tau_1| Y_{\tau_1}^\varepsilon\in U_1\cup \partial U_1)}
{\ln(\varepsilon^{-1})}\leq \dfrac{1}{\gamma_1} \ .
\end{equation}

The above \eqref{Thm:MajorResultSGDPerturbationConvergenceTime:Proof:SigmatildeLargerSigma} and \eqref{Thm:MajorResultSGDPerturbationConvergenceTime:Proof:SigmatildeMinusTau} together indicate that

\begin{equation}\label{Thm:MajorResultSGDPerturbationConvergenceTime:Proof:SigmaMinusTauPart1}
\lim\limits_{\varepsilon\rightarrow 0}
\dfrac{\mathbf{E}_x (\sigma_1-\tau_1| Y_{\tau_1}^\varepsilon\in U_1\cup \partial U_1)}
{\ln(\varepsilon^{-1})}\leq \dfrac{1}{\gamma_1} \ .
\end{equation}

Let us then turn to the case when $Y_{\tau_1}^\varepsilon \not\in U_1\cup \partial U_1$. Notice that by definition of $\tau_1$, in this case we have $Y_{\tau_1}^\varepsilon\in C_1(h/2)$. Thus by part (iv) of Lemma \ref{Lm:GeometricLemma}, we can construct a tubular neighborhood containing $Y_{\tau_1}^\varepsilon$ and connecting $C_1(h/2)$ to $C_1(-h)$. The neighborhood can be chosen so that its boundaries are of a positive distance to $Y_{\tau_1}^\varepsilon$. Hence we can then apply Lemma \ref{Lm:FiniteTimeComparisonGradientFlowANDSGDPerturbation} again so that
for any $\nu>0$, there exist some $\varepsilon_0>0$, for any $0<\varepsilon<\varepsilon_0$ we have

\begin{equation}\label{Thm:MajorResultSGDPerturbationConvergenceTime:Proof:SigmaMinusTauPart2}
\dfrac{\sup\limits_{H-r>F(y)\geq F(O_1)+h/2}\mathbf{E}_y (\sigma_1-\tau_1| Y_{\tau_1}^\varepsilon \not\in  U_1\cup \partial U_1)}
{\ln(\varepsilon^{-1})}
\leq  \nu \ .
\end{equation}

Combining equations \eqref{Thm:MajorResultSGDPerturbationConvergenceTime:Proof:SigmaMinusTauPart1}, \eqref{Thm:MajorResultSGDPerturbationConvergenceTime:Proof:SigmaMinusTauPart2} we see that, for any $\nu>0$, there
exist some $\varepsilon_0>0$ so that for any $0<\varepsilon<\varepsilon_0$ we have

\begin{equation}\label{Eq:Sigma2-Tau1ConditionalCombined}
\begin{array}{ll}
&\dfrac{\sup\limits_{H-r>F(y)\geq F(O_1)+h/2}\mathbf{E}_y(\sigma_1-\tau_1)}{\ln(\varepsilon^{-1})}
\\
= & \dfrac{\sup\limits_{H-r>F(y)\geq F(O_1)+h/2}\mathbf{E}_y (\sigma_1-\tau_1| Y_{\tau_1}^\varepsilon \not\in  U_1\cup \partial U_1)
\cdot \mathbf{P}_y(Y_{\tau_1}^\varepsilon \not\in  U_1\cup \partial U_1)}
{\ln(\varepsilon^{-1})}
\\
& + \dfrac{\sup\limits_{H-r>F(y)\geq F(O_1)+h/2}
\mathbf{E}_y (\sigma_1-\tau_1| Y_{\tau_1}^\varepsilon\in U_1\cup \partial U_1)
\cdot \mathbf{P}_y( Y_{\tau_1}^\varepsilon\in U_1\cup \partial U_1)}
{\ln(\varepsilon^{-1})}
\\
\leq &  \dfrac{1}{\gamma_1}+\nu \ .
\end{array}
\end{equation}

Since $\nu>0$ can be picked arbitrarily small as $\varepsilon>0$ is small, we see from \eqref{Eq:Sigma2-Tau1ConditionalCombined}
that we finally have

\begin{equation}\label{Eq:Sigma2-Tau1FinalEstimate}
\lim\limits_{\varepsilon\rightarrow 0}\dfrac{\sup\limits_{H-r>F(y)\geq F(O_1)+h/2}\mathbf{E}_y(\sigma_1-\tau_1)}{\ln(\varepsilon^{-1})}
\leq \dfrac{1}{\gamma_1} \ .
\end{equation}

\

Notice that we have $F(Y^\varepsilon_{\sigma_1})\leq F(O_1)-h$.
Due to strong Markov property, given $Y_{\sigma_1}^\varepsilon$, the process
$\widetilde{Y}_s^\varepsilon=Y_{\sigma_1+s}^\varepsilon$ is an independent copy of $Y_t^{\varepsilon}$. Therefore we can iteratively make use of
 \eqref{Eq:Sigma2-Tau1FinalEstimate}. From Lemma \ref{Lm:FiniteTimeComparisonGradientFlowANDSGDPerturbation}
 we know that, given any small $\rho>0$, as $\varepsilon>0$ is small, with probability greater or equal than
 $1-\rho$, $Y_{\tau_2}^\varepsilon$ lie on $C_2(h/2)$, i.e.,
$F(Y_{\tau_2}^\varepsilon)=F(O_2)+h/2$. Similarly, given any small $\rho>0$,
as $\varepsilon>0$ is small, with probability greater or equal than
 $1-\rho$, $Y_{\sigma_2}^\varepsilon$ lie on $C_2(-h)$, i.e., $F(Y_{\sigma_2}^\varepsilon)=F(O_2)-h$.

Thus by running the above argument iteratively, we see that for any $\rho>0$ sufficiently small,
there exist some $\varepsilon_0>0$ such that for any $0<\varepsilon<\varepsilon_0$, with probability greater or equal than $1-\rho$,
the random variable $Y_{\sigma_k}^\varepsilon$ satisfy $F(Y_{\sigma_k}^\varepsilon)=F(O_k)-h$. Conditioned on this event, we have
$$\lim\limits_{\varepsilon\rightarrow 0}\dfrac{\mathbf{E}_x T_x^\varepsilon}{\ln(\varepsilon^{-1})}
\leq \dfrac{k}{\gamma_1} \ .$$
So we arrive at the conclusion of this Theorem. $\square$

\

Taking into account that $Y_t^{\sqrt{\beta}}=X_{t/\beta}$ for the process $X_t$ defined in \eqref{Eq:SGD} and the process $Y_t^\varepsilon$ defined in \eqref{Eq:SGDPerturbation}, we conclude that
$\tau^\beta$
in \eqref{Eq:HittingTimeSGDManySaddle} satisfies $\tau^\beta=T_x^{\sqrt{\beta}}/\beta$. Thus
we see that we have the following corollary in regards to the diffusion approximation of the
stochastic gradient descent (Section 1.3).

\begin{corollary}\label{Corollary:MajorResultSGDConvergenceTime}
Consider the process $X_t$ defined as in \eqref{Eq:SGD} and let $F(\bullet)$ satisfy the landscape conditions as in Theorem \ref{Thm:MajorResultSGDPerturbationConvergenceTime}, and also suppose the noise in \eqref{Eq:SGD} is sufficiently omnidirectional in the sense that $\sigma(x)\sigma^T(x) $ is uniformly positive definite. Let $x^*$ be the unique local minimum of $F$ as in Theorem \ref{Thm:MajorResultSGDPerturbationConvergenceTime}, then

(i) For any small $\rho>0$, with probability at least $1-\rho$, SGD diffusion process $X_t$ in \eqref{Eq:SGD} converges
to the minimizer $x^*$ for sufficiently small $\beta$ after passing through all $k$ saddle points $O_1$, ..., $O_k$;

(ii) Consider the stopping time $\tau^\beta$ defined in
\eqref{Eq:HittingTimeSGDManySaddle}. Then as $\beta\downarrow 0$, conditioned on the above convergence of $SGD$ diffusion process $X_t$,
we have
\begin{equation}\label{Eq:TotalExitTimeAsymptoticSGDManySaddle}
\lim\limits_{\beta \rightarrow 0}\dfrac{\mathbf{E} \tau^\beta}{\beta^{-1}\ln \beta^{-1}}\leq
\dfrac{k}{2\gamma_1} \ .
\end{equation}
\end{corollary}

\section*{Acknowledgement} We thank the anonymous referee for valuable comments that help to improve a preliminary version of this paper.

\bibliographystyle{plain}
\bibliography{bibliography_SGDescape}

\begin{thebibliography}{10}

\bibitem{[Arnold]}
V.I. Arnold.
\newblock Geometric {M}ethods in the {T}heory of {O}rdinary {D}ifferential
  {E}quations.
\newblock {\em Grundlehren der mathematischen Wissenschaften}, 250, 1988.

\bibitem{Bach-MoulinesSGD2013NIPS}
F.~Bach and E.~Moulines.
\newblock Non--strongly--convex smooth stochastic approximation with
  convergence rate $o(1/n)$.
\newblock {\em Neural Information and Processing Systems}, pages 773--781,
  2013.

\bibitem{[Bakhtin2008SPA]}
Y.~Bakhtin.
\newblock Exit asymptotics for small diffusion about an unstable equilibrium.
\newblock {\em Stochastic Processes and their Applications}, 118:839--851,
  2008.

\bibitem{[BakhtinHeteroclinic]}
Y.~Bakhtin.
\newblock Noisy heteroclinic networks.
\newblock {\em Probability Theory and Related Fields}, 150(1-2):1--42, 2011.

\bibitem{[StrengthenedHartmanGrobman]}
G.~Belitskii and V.~Rayskin.
\newblock On the {G}robmann--{H}artman {T}heorem in $\alpha$--{H}\"{o}lder
  class for {B}anach spaces.
\newblock {\em \verb"https://www.ma.utexas.edu/mp_arc/c/11/11-134.pdf"}.

\bibitem{BenvenisteEtAlAdaptiveAlgorithms}
A.~Benveniste, M.~M\'{e}tivier, and P.~Priouret.
\newblock {\em Adaptive Algorithms and Stochastic Approximations, Applications
  of Mathematics, 22}.
\newblock Springer, 1990.

\bibitem{[BORKAR]}
Vivek~S. Borkar.
\newblock {\em Stochastic Approximation: A Dynamical Systems Viewpoint}.
\newblock Cambridge University Press, 2008.

\bibitem{BottouSGD2010}
L.~Bottou.
\newblock Large-scale machine learning with stochastic gradient descent.
\newblock {\em Proceedings of COMPSTAT'2010}, pages 177--186, 2010.

\bibitem{[Brin-Stuck]}
M.~Brin and G.C. Stuck.
\newblock {\em Introduction to dynamical systems}.
\newblock Cambridge University Press, 2002.

\bibitem{Flammarion-BachCOLT2015}
N.~Flammarion and F.~Bach.
\newblock From averaging to accelaration, there is only a step--size.
\newblock {\em Conference on Learning Theory (COLT)}, 2015.

\bibitem{[FWbook]}
M.~Freidlin and A.~Wentzell.
\newblock {\em Random perturbations of dynamical systems, Second Edition}.
\newblock Springer, 1998.

\bibitem{[GeEtAl1503]}
R.~Ge, F.~Huang, C.~Jin, and Y.~Yuan.
\newblock scaping from saddle points--online stochastic gradient descent for
  tensor decomposition.
\newblock {\em Proceedings of The 28th Conference on Learning Theory (COLT)},
  pages 797--842, 2015.

\bibitem{[HartmanODE]}
P.~Hartman.
\newblock {\em Ordinary Differential Equations}.
\newblock John Wiley and Sons, New York, 1964.

\bibitem{JunchiEtAlSDE2017}
W.~Hu, C.J. Li, L.~Li, and J.G. Liu.
\newblock On the diffusion approximation of nonconvex stochastic gradient
  descent.
\newblock {\em Annals of Mathematical Science and Applications}, 4(1):3--32,
  2019.

\bibitem{[KatokHasselblatt]}
A.~Katok and B.~Hasselblatt.
\newblock Introduction to the modern theory of dynamical systems.
\newblock {\em Encyclopaedia of Mathematics and its Applications}, 54, 1995.

\bibitem{[Kifer1981]}
Y.~Kifer.
\newblock The exit problem for small random perturbations of dynamical systems
  with a hyperbolic fixed point.
\newblock {\em Israel Journal of Mathematics}, 40(1):74--96.

\bibitem{[Koralov2004]}
L.~Koralov.
\newblock Random perturbations of 2--dimensional hamiltonian flows.
\newblock {\em Probability Theory and Related Fields}, 129:37--62, 2004.

\bibitem{Kushner-Yin}
H.J. Kushner and G.~George~Yin.
\newblock {\em Stochastic approximation and recursive algorithms and
  applications, Stochastic Modeling and Applied Probability, 35}.
\newblock Springer, 2003.

\bibitem{lee2016gradient}
J.D. Lee, M.~Simchowitz, M.I. Jordan, and B.~Recht.
\newblock Gradient descent only converges to minimizers.
\newblock {\em Proceedings of The 29th Conference on Learning Theory (COLT)},
  pages 1246--1257, 2016.

\bibitem{WeinanEtAlSDE2017}
Q.~Li, C.~Tai, and W.~E.
\newblock Stochastic modified equations and adaptive stochastic gradient
  algorithms.
\newblock {\em International Conference on Machine Learning (ICML)}, pages
  2101--2110, 2017.

\bibitem{WeinanEtAlSDE2018}
Q.~Li, C.~Tai, and W.~E.
\newblock Stochastic {M}odified {E}quations and {D}ynamics of {S}tochastic
  {G}radient {A}lgorithms {I}: {M}athematical {F}oundations.
\newblock {\em arXiv:1811.01558[cs.LG]}, 2018.

\bibitem{[Milnor1963]}
J.~Milnor.
\newblock {\em Morse Lemma}.
\newblock Princeton University Press, 1963.

\bibitem{[MonterBakhtinNormalForm]}
S.A.A. Monter and Y.~Bakhtin.
\newblock Normal forms approach to diffusion near hyperbolic equilibira.
\newblock {\em Nonlinearity}, 24:1883--1907, 2011.

\bibitem{NonGaussianitySGD-Noise}
A.~Panigrahi, R.~Somani, N.~Goyal, and P.~Netrapalli.
\newblock Non-{G}aussianity of {S}tochastic {G}radient {N}oise.
\newblock {\em Science meets Engineering of Deep Learning (SEDL) workshop, 33rd
  Conference on Neural Information Processing Systems (NeurIPS 2019),
  Vancouver, Canada}, 2019.

\bibitem{pemantle1990nonconvergence}
R.~Pemantle.
\newblock Nonconvergence to unstable points in urn models and stochastic
  approximations.
\newblock {\em The Annals of Probability}, 18(2):698--712, 1990.

\bibitem{Robbins-Monroe}
H.~Robbins and S.~Monro.
\newblock A {S}tochastic {A}pproximation {M}ethod.
\newblock {\em The Annals of Mathematical Statistics}, 22(3), 1951.

\bibitem{simsekli2019tail}
U.~Simsekli, L.~Sagun, and M.~Gurbuzbalaban.
\newblock A tail-index analysis of stochastic gradient noise in deep neural
  networks.
\newblock {\em International Conference on Machine Learning (ICML)}, 2019.

\bibitem{[SGD-continuous-time]}
J.~Sirignano and K.~Spiliopoulos.
\newblock Stochastic gradient descent in continuous time.
\newblock {\em SIAM Journal on Financial Mathematics}, 8(1):933--961, 2017.

\bibitem{[sun2015nonconvex]}
J.~Sun, Q.~Qu, and J.~Wright.
\newblock When are nonconvex problems not scary?
\newblock {\em Preprint, arXiv:1510.06096[math.OC]}.

\bibitem{[Pathway-Energy-Landscape]}
J.~Yin, Y.~Wang, Jeff~Z.Y. Chen, P.~Zhang, and L.~Zhang.
\newblock Construction of a {P}athway {M}ap on a {C}omplicated {E}nergy
  {L}andscape.
\newblock {\em Phys. Rev. Lett.}, 124(090601), 2, March 2020.

\bibitem{ZhanxingSGDAnisotropic}
Z.~Zhu, J.~Wu, B.~Yu, B.~Wu, and J.~Ma.
\newblock The anisotropic noise in stochastic gradient descent: Its behavior of
  escaping from minima and regularization effects.
\newblock {\em International Conference on Machine Learning (ICML) 2019,
  arXiv:1803.00195[stat.ML]}, 2019.

\end{thebibliography}

\end{document}